# Efficient multi-level *hp*-finite elements in arbitrary dimensions

Philipp Kopp[a,*], Ernst Rank[b,c], Victor M. Calo[d], Stefan Kollmannsberger[a]

[a] *Chair of Computational Modeling and Simulation, Technical University of Munich, Germany*
[b] *Chair for Computation in Engineering, Technical University of Munich, Germany*
[c] *Institute for Advanced Study, Technical University of Munich, Germany*
[d] *Applied Mathematics, School of Electrical Engineering, Computing & Mathematical Science, Faculty of Science and Engineering, Curtin University, Australia*

**Abstract**

We present an efficient algorithmic framework for constructing multi-level *hp*-bases that uses a data-oriented approach that easily extends to any number of dimensions and provides a natural framework for performance-optimized implementations. We only operate on the bounding faces of finite elements without considering their lower-dimensional topological features and demonstrate the potential of the presented methods using a newly written open-source library. First, we analyze a Fichera corner and show that the framework does not increase runtime and memory consumption when compared against the classical *p*-version of the finite element method. Then, we compute a transient example with dynamic refinement and derefinement, where we also obtain the expected convergence rates and excellent performance in computing time and memory usage.

## 1. Introduction

For smooth solutions, *p*-finite element methods converge exponentially with respect to the number of unknowns [1]; for rough solutions, *hp*-finite elements recover the exponential convergence [2,3] where element sizes and polynomial degrees gradually decrease towards the solution singularities. The decision of whether to refine $h$ or $p$ can be automated by combining an error estimator and a smoothness indicator (see, e.g., [4,5]). Classical implementations of *hp*-finite element methods follow a refinement by replacement strategy [6,7], in which smaller elements substitute elements marked for refinement. Once the new mesh is constructed, the shape functions of the smaller elements must be connected to the shape functions of neighboring elements in a non-trivial way to obtain a globally $C^0$ continuous finite element basis. In particular, the richer space of the fine elements must be constrained to the coarse shape functions of the neighbors. Especially in regions with sharp transitions between coarse and fine elements, these constraints introduce complicated dependencies on already constrained elements.

We follow the multi-level *hp*-method introduced in [8–11] that simplifies the construction of *hp*-bases using a refine-by-superposition approach rather than a refinement-by-replacement. This superposition defines shape functions on all levels of the refinement tree and simplifies the construction of compatible basis functions ($C^0$ continuous across element interfaces) significantly, at the cost of slightly larger shape function supports in transition zones. The original framework [8,9] builds an object-oriented data structure to store the various topological relations between the elements and their sub-topologies (faces, edges, nodes in three dimensions). Then, basis functions are associated to topological components, which allows the formulation of two rules for constructing multi-level *hp*-bases.

* Corresponding author.
*E-mail address:* philipp.kopp@tum.de (P. Kopp).

Motivated by the simplicity of the method, we extend the framework to four dimensions to allow conforming space–time discretizations of transient three-dimensional problems. The application in space–time discretizations of parabolic problems follows in [12]. Here, we focus on extending the multi-level *hp*-method to higher dimensions. While being convenient, the original object-oriented data structure becomes increasingly complex with each new dimension. The classical multi-level *hp*-approach in four dimensions constructs objects for each 4-cube and its eight cubes, twenty-four faces, thirty-two edges, and sixteen nodes. Each topological entity lists its adjacent finite elements (4-cubes) and its overlay topologies (e.g., faces may list four overlay faces) and its topological elements (e.g., faces list their four edges and four nodes, while edge list their two nodes).

The number of dependencies in such a data structure consumes a significant amount of memory while being tedious to construct and maintain. Herein, we introduce a data-oriented alternative. While the basis functions are identical, we only require basic topological relations between cells in the refinement tree, not their lower-dimensional bounding topologies. Using these basic relations, we introduce a four-step algorithm to construct a $d$-dimensional mask of Boolean values; this tensor-product mask activates the proper shape functions from the entire tensor-product space of integrated Legendre shape functions leaving the others inactive. We construct location maps that combine the active shape functions on each element to the globally compatible *hp*-basis functions. By storing all the data in array-type containers, the effort of implementing a sophisticated object-oriented data structure becomes developing a capable low-level implementation of these containers. NumPy has been a significant influence on the ideas presented here. Moreover, we want to mention [13,14] as approaches that relate closely to our work. Because adaptive algorithms based on a posteriori error estimators have already been investigated with the multi-level *hp*-method [4,5], we focus on our simple, robust, and dimension independent multi-level *hp*-algorithms and data structures. Our method yields almost identical basis functions compared to [8,9], with the exception of slightly smaller supports of some linear modes.

Throughout this paper, we use the linear heat equation as a model problem:

$$\begin{aligned} -\nabla \cdot (\kappa \nabla u) &= f && \text{on } \Omega \\ u &= u_d && \text{on } \Gamma_D \\ \kappa \nabla u \cdot n &= h && \text{on } \Gamma_N, \end{aligned} \tag{1}$$

where $\Gamma_D \cap \Gamma_N = \partial\Omega$ and $\Gamma_D \cup \Gamma_N = \emptyset$. The weak form of (1) reads as follows: Find $u \in u_d + \mathcal{H}_0^1(\Omega)$, such that

$$\int_\Omega \nabla w \cdot \kappa \nabla u \, \mathrm{d}\Omega = \int_\Omega wf \, \mathrm{d}\Omega + \int_{\Gamma_N} wh \, \mathrm{d}\Gamma_N \qquad \forall w \in \mathcal{H}_0^1(\Omega).$$

We assume that we can mesh $\Omega$ using a Cartesian grid, which is possible in our application of the presented method to space–time finite elements in additive manufacturing. This assumption simplifies some of the algorithms but is not a restriction of the method as the ideas can be generalized to unstructured meshes. Moreover, despite significant effort in the past [15–17], mesh generation for higher-order methods is still challenging and being an area of ongoing research, which applies even more to higher dimensions. In particular, automatic *hp*-refinements require sufficiently accurate geometric approximations when refining elements on the boundary, which complicates the use of mesh generators that use polynomial approximations. A common alternative for complicated cases is to use immersed techniques [1,18,19] that use non-conforming (Cartesian) grids. The influence from outside of the physical domain is then minimized in the numerical integration of the element stiffness matrices, such that accuracies comparable to boundary conforming discretizations can be achieved.

This paper is structured as follows. Section 2 introduces our new algorithms and data structures for single *p*-finite element meshes without refinement. We first discuss tensor-products of integrated Legendre shape functions on each element that we can selectively activate or deactivate using *tensor-product masks*. We define these masks on each element and construct them in Section 2.2 such that the corresponding active shape functions can be combined to globally compatible basis functions in Section 2.3. We then introduce hierarchical space-tree-based overlays in Section 3 and extend our algorithms for constructing the multi-level *hp* tensor-product masks in Section 3.2 and for combining the active shape functions to global basis functions in Section 3.3. Finally, we verify the algorithms and their implementation on a set of examples in Section 4 and show that they behave optimally.

We provide an open-source licenced reference implementation in Python for all examples (including Sections 2 and 3) at https://gitlab.com/hpfem/publications/2021_nd-mlhp. We assess the performance potential of our method by computing our examples using an optimized implementation in C++ that is also available in the same repository.

## 2. Single-level *p*-finite element basis

In this section, we consider *p*-finite element meshes without refinements and with varying polynomial degrees per element and coordinate direction. This case is well known by the community and allows us to introduce our approach in a familiar setting. Our meshes consist of *d*-dimensional (hyper-) cuboids, i.e., line segments in one dimension, quadrilaterals in two dimensions, cuboids or bricks in three dimensions and so on. We assign each element a unique id between 0 and $n_{\mathrm{el}} - 1$ in an arbitrary order. Locally, all elements have two faces per coordinate direction, some of which are interfaces to other elements. We store the ids of both neighbors in $d$ directions in an array $N$ of size $(n_{el}, d, 2)$, where we use a dedicated *no cell* value for domain boundary faces; for example, $-1$ for signed integers or the maximum representable number for unsigned integers. Unlike conventional *p*- and *hp*-finite element approaches, we do not store our mesh as graph of vertices, edges, faces, and so on; instead, we formulate our algorithms such that they only use basic topological relations between elements.

The main challenge in constructing $C^0$ continuous basis functions is to handle neighboring elements with varying polynomial degrees. To simplify the description of our algorithms, we assume that the coordinate axes of neighboring elements are aligned, which is generally not true for unstructured meshes.

### 2.1. Integrated Legendre shape functions

We begin by discussing the tensor-products of integrated Legendre polynomials in the local coordinates $r_0, \dots, r_{d-1}$, with $r_i \in [-1, 1]$ and their selective activation or deactivation. First, we define the integrated Legendre polynomials of degree $q$ in one dimension

$$I_0(r) = \frac{1}{2}(1-r) \qquad I_1(r) = \frac{1}{2}(1+r) \qquad I_q(r) = \frac{L_q(r) - L_{q-2}(r)}{\sqrt{4q-2}} \quad \text{for } q > 1, \tag{2a}$$

$$I_0'(r) = -\frac{1}{2} \qquad I_1'(r) = \frac{1}{2} \qquad I_q'(r) = \frac{L_q'(r) - L_{q-2}'(r)}{\sqrt{4q-2}} \quad \text{for } q > 1, \tag{2b}$$

using the Legendre polynomials

$$L_0(r) = 1 \qquad L_1(r) = r \qquad L_q(r) = \frac{2q-1}{q} r\, L_{q-1} - \frac{q-1}{q} L_{q-2} \quad \text{for } q > 1, \tag{3a}$$

$$L_0'(r) = 0 \qquad L_1'(r) = 1 \qquad L_q'(r) = \frac{2q-1}{q}\left(L_{q-1} + r L_{q-1}'\right) - \frac{q-1}{q} L_{q-2}' \quad \text{for } q > 1. \tag{3b}$$

The hierarchical nature of this basis yields the important property that the basis $\{I_0, \dots, I_p\}$ for $\mathcal{P}_p([-1, 1])$ can be extended to form a basis for $\mathcal{P}_{p+1}([-1, 1])$ by simply adding the function $I_{p+1}$. Another important property is that $I_0$ is the only non-zero function on the left boundary. Similarly, $I_1$ is the only non-zero function on the right boundary, which allows us to classify the functions into left and right nodal modes ($I_0$ and $I_1$) and internal or bubble modes ($I_p$, if $p \geq 2$). Fig. 1 shows the integrated Legendre bases for $p = 1 - 3$.

In higher dimensions, we define an integrated Legendre basis in each coordinate direction: $\{I_0, \dots, I_{p_i}\}$. The polynomial degrees in different coordinate axes may vary; thus, $p$ is now a tuple $(p_0, \dots, p_{d-1})$ with a possibly different polynomial degree per axis. This feature, for example, allows choosing different polynomial degrees in space and time for finite element discretizations of transient problems, or it can reduce the number of unknowns in combination with a capable error estimator. Now, a basis for the *tensor-product space* can be defined as follows:

$$\left\{ N_\alpha(r) = I_{\alpha_0}(r_0) \cdot I_{\alpha_1}(r_1) \cdot \ \dots \ \cdot I_{\alpha_{d-1}}(r_{d-1}) \;\middle|\; \text{for } \alpha \in \{0, \dots, p_0\} \times \dots \times \{0, \dots, p_{d-1}\} \right\}, \tag{4}$$

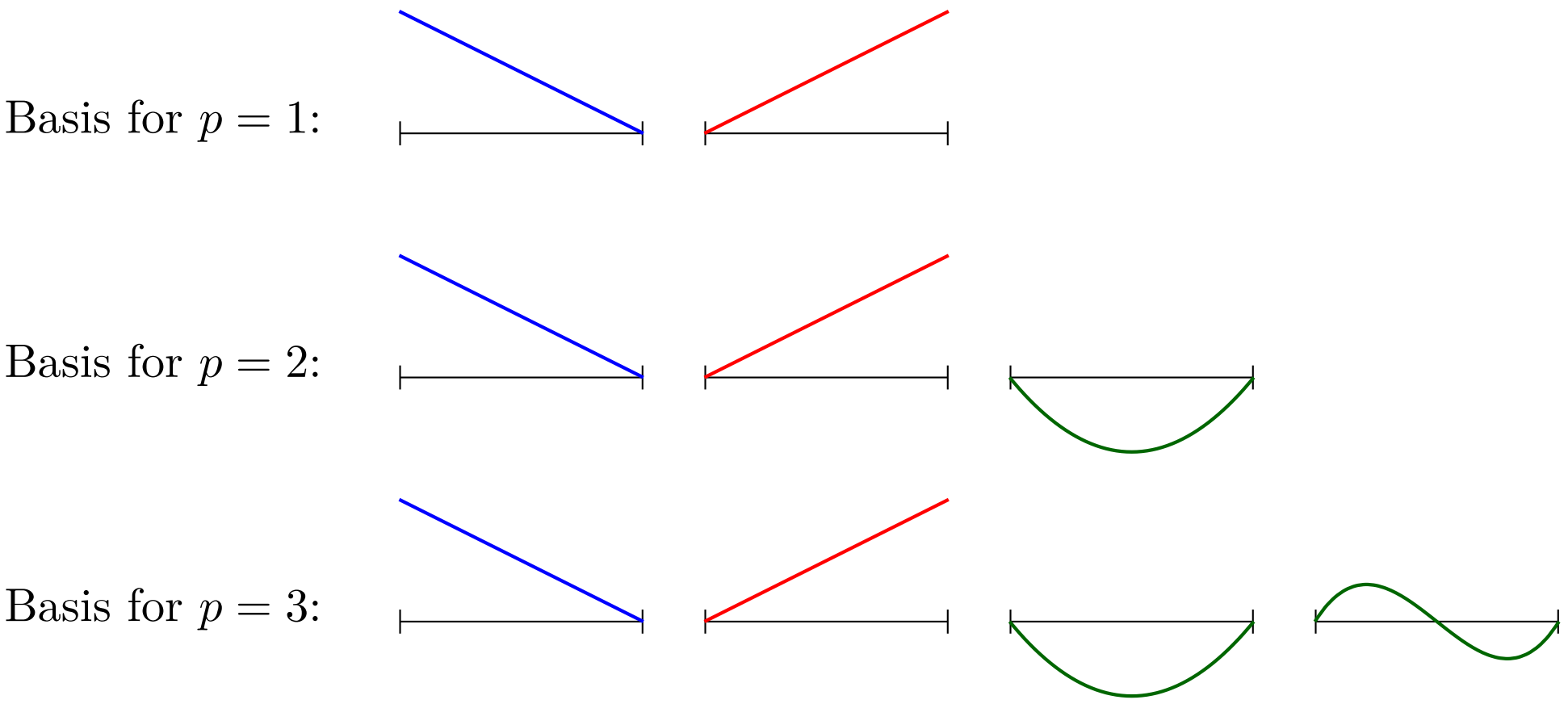

**Fig. 1.** One-dimensional Integrated Legendre shape functions.

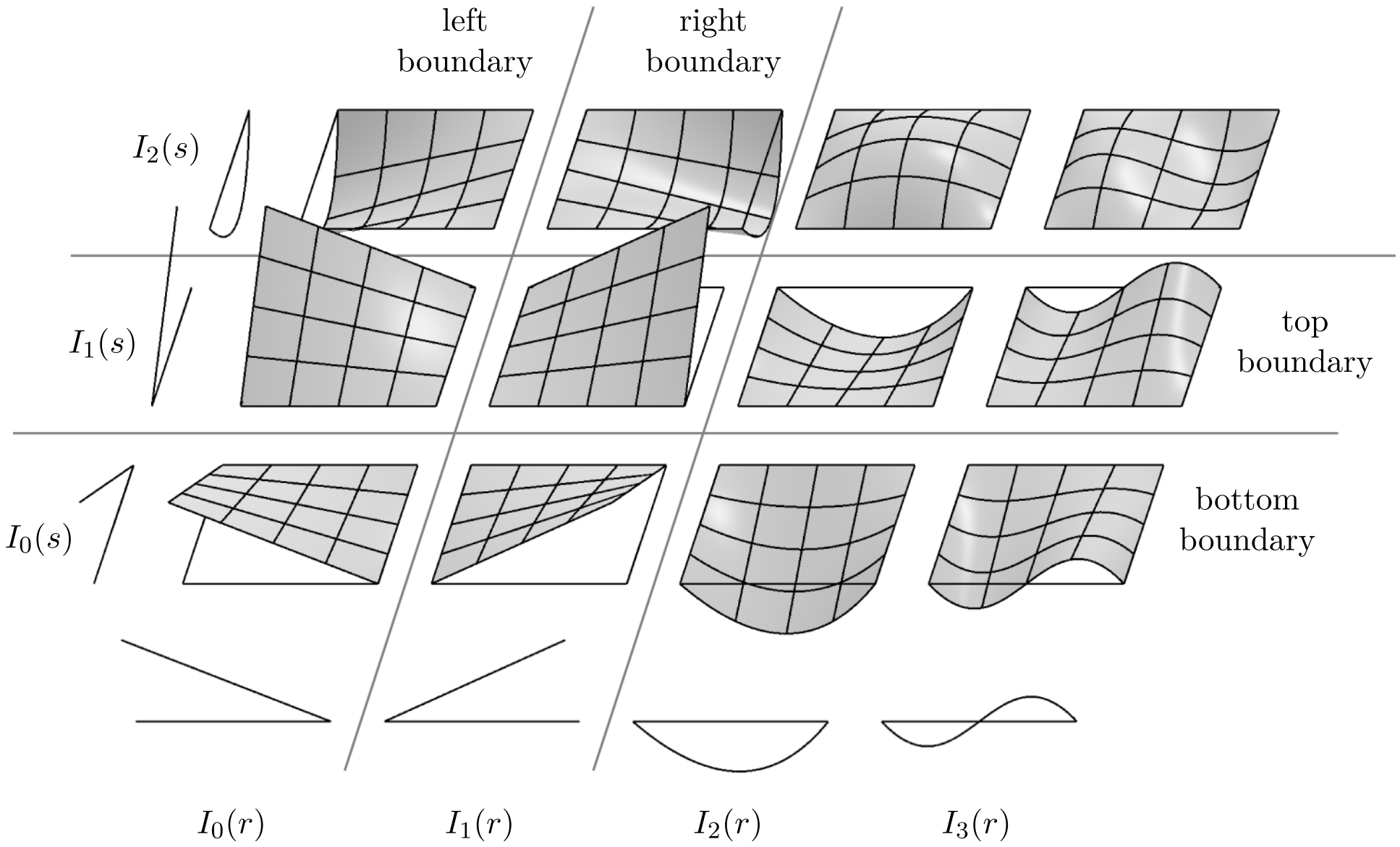

**Fig. 2.** Tensor-product of Integrated Legendre shape functions for $p = (3, 2)$.

where $r_i$ are the local coordinates in direction $i$, and each $\alpha_i \in \{0, \ldots, p_i\}$ indexes the one-dimensional shape functions in direction $i$. In two dimensions, for example, (4) becomes

$$\left\{ N_{ij}(r_0, r_1) = I_i(r_0) \cdot I_j(r_1) \;\middle|\; \text{for } i \in \{0, \ldots, p_0\},\ j \in \{0, \ldots, p_1\} \right\}, \tag{5}$$

where $r_0$ and $r_1$ are the two local coordinates. The functions in (4) and (5) are defined locally on each element and they span complete polynomial spaces in the local coordinate systems. From now on, we refer to them as shape functions to distinguish them from the finite element basis functions that we construct later by "gluing" together several shape functions from adjacent elements.

Fig. 2 displays the two-dimensional shape functions in (5) for $p = (3, 2)$. Independent of the polynomial degree, we identify the first and second columns for functions with support on the left and right boundaries and the first and second rows for functions with support on the bottom and top boundaries. More generally, we associate a $d - 1$ dimensional *slice* of the tensor-product to each boundary of a cell. In three dimensions, we have a two-dimensional slice for each of the six boundaries; in four dimensions, we have a three-dimensional slice for each of the eight boundaries, and so on. Together with the hierarchy in $p$, this structure forms the basis of our algorithms for

$$\begin{bmatrix} 1 & 1 & 1 & 1 \\ 1 & 1 & 1 & 1 \\ 1 & 1 & 1 & 1 \end{bmatrix}$$

(a) All functions active (Figure 2).

$$\begin{bmatrix} 1 & 0 & 1 & 1 \\ 1 & 0 & 1 & 1 \\ 1 & 0 & 1 & 1 \end{bmatrix}$$

(b) Right side deactivated.

$$\begin{bmatrix} 1 & 0 & 1 & 1 \\ 1 & 0 & 1 & 1 \\ 0 & 0 & 0 & 0 \end{bmatrix}$$

(c) Bottom and right sides deactivated.

$$\begin{bmatrix} 1 & 0 & 0 & 0 \\ 1 & 1 & 1 & 1 \\ 1 & 0 & 0 & 0 \end{bmatrix}$$

(d) Only top and left sides activated.

**Fig. 3.** Tensor-product masks in two dimensions.

connecting the shape functions into $p$- and $hp$-basis functions. For example, in Section 3.2, we impose homogeneous Dirichlet conditions on the boundaries of overlay cells by simply deactivating the corresponding slices in the tensor-product. We now formalize our notation by identifying slices in the tensor-product by their normal axis and the index along this axis (either 0 or 1). For example, in two dimensions, $(0, 1)$ corresponds to the right boundary slice, which is normal to axis zero at index one, and $(1, 0)$ corresponds to the bottom boundary, which is normal to axis one at index zero. Similarly, the slice $(0, 0)$ and $(1, 1)$ correspond to the left and top boundaries. In three dimensions we have $(0, 0)$, $(0, 1)$, $(1, 0)$, $(1, 1)$, $(2, 0)$, and $(2, 1)$ for the functions with support on the left, right, front, back, bottom, and top faces, respectively. Our algorithms extend well to higher dimensions because we do not introduce a further classification of the shape functions into nodal modes, edge modes, face modes, and so on, unlike traditional $p$- and $hp$-finite element approaches.

The first step in creating our global $p$- and $hp$-basis functions is to select the right shape functions from the tensor-product of each element, such that we can later glue them together without gaps. To keep track of active and inactive shape functions in the tensor-product, we introduce a $d$-dimensional tensor of Boolean values for each element that we call a *tensor-product mask*. An entry $M_\alpha$, with $\alpha = (\alpha_0, \ldots, \alpha_{d-1})$ in this tensor-product mask determines whether the shape function $N_\alpha$ is active or not. All shape functions with indices outside the shape of $M$ are implicitly inactive. Therefore, $M$ must be large enough to include all active shape functions. Fig. 3(a) shows the tensor-product mask corresponding to Fig. 2, where all shape functions are active. We then impose homogeneous Dirichlet boundary conditions on the right side in Fig. 3(b) and then additionally on the bottom side in Fig. 3(c). Fig. 3(d) shows the result of activating the left and top boundary functions starting from a cell without any active functions.

**Remark 2.1.** The integrated Legendre polynomials (2) can be expanded and evaluated using Horner's scheme. However, the recursive definition we use is computationally more efficient when we compute all shape functions $I_0 - I_p$ at once.

## *2.2. Compatible tensor-product masks*

We can assign individual polynomial degrees per element and direction to our $p$-finite element bases, which complicates creating compatible basis functions from connected elements with different polynomial degrees in the interface directions. There are generally two options for resolving shape functions that exist on one side of an element interface that are left without counterpart on the other side. The *minimum degree strategy* deactivates shape functions without counterpart across the interface and the *maximum degree strategy* activates functions on the other side instead. In terms of our tensor-product masks, the corresponding slices of both elements are initially incompatible. We restore compatibility by comparing the values in both slices using a component-wise logical *and* operation or alternatively a component-wise logical *or* operation. We use the result of this comparison to overwrite the original Boolean values on both sides. Fig. 4 shows several examples where we restore compatibility using logical operations on the interface slices. The logical *and* results in a minimum degree while the logical *or* results in a maximum degree strategy, but they can be used more generally to make tensor-product masks compatible, whether they originate from different polynomial degrees or not (see Section 3.2).

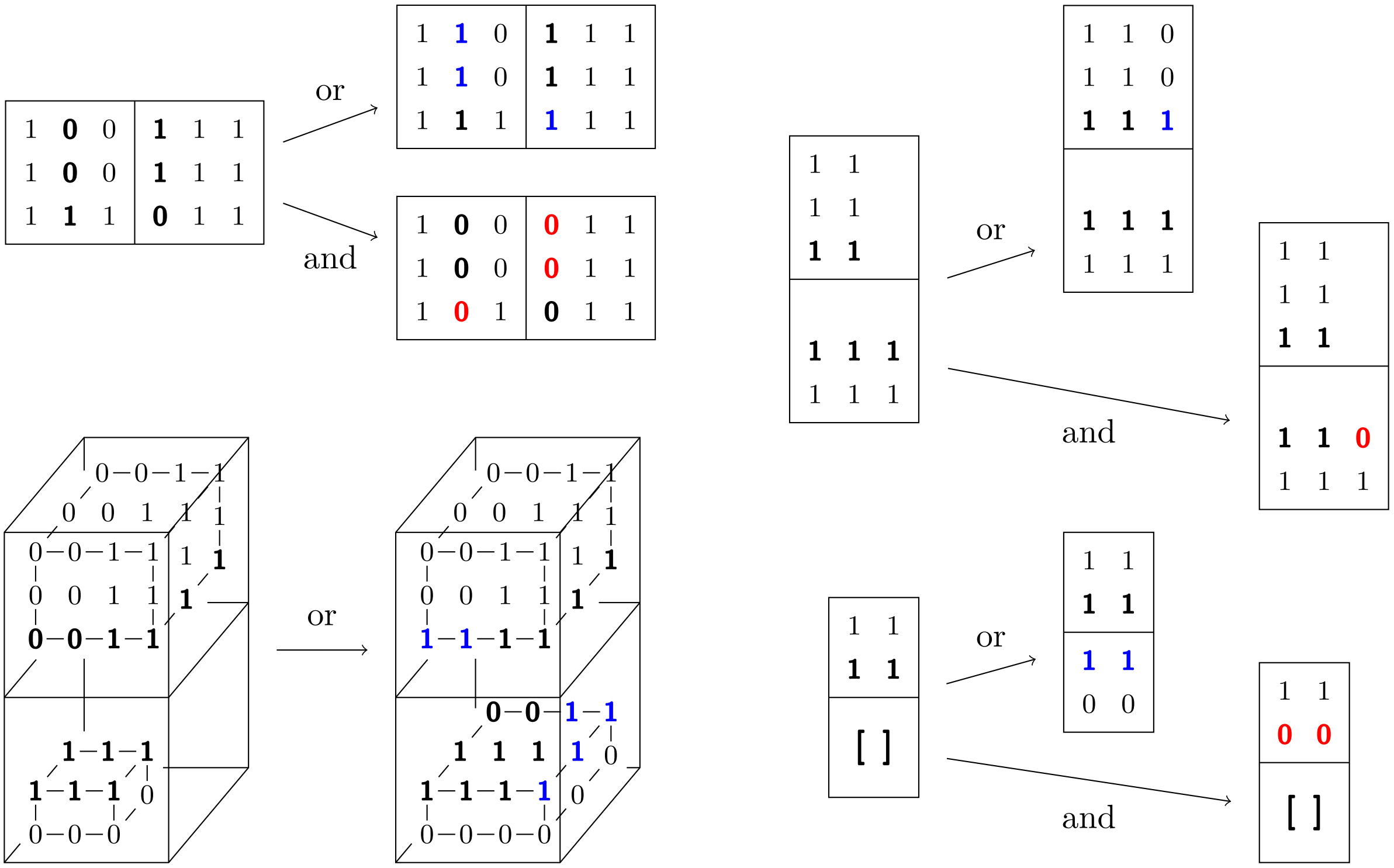

**Fig. 4.** Restore interface compatibility between the tensor-product masks of two cells by comparing and overwriting the interfaces slices using logical operations. Slices to be compared are printed in boldface; blue entries were activated due to logical *or*; red entries were deactivated due to logical *and*. (For interpretation of the references to color in this figure legend, the reader is referred to the web version of this article.)

In this setting, we construct the tensor-product masks for $p$-finite element meshes in two steps. First, we initialize all masks by activating the full tensor-product up to the polynomial degrees assigned to the elements. Then, we force compatibility on all element interfaces using logical *and* operations if we want minimum interface degrees or logical *or* operations if we want maximum degrees. Fig. 5 demonstrates this procedure on a two-dimensional mesh. Despite being compatible, some pairs of tensor-product masks have different sizes in the interface direction, which uses the definition that entries are inactive if they are not explicitly stored.

In $d$ dimensions, we execute the second step $d-1$ times to reach all elements that our shape functions connect to. For example, in three dimensions the basis functions associated to an edge may be connected from four elements; hence, for information to travel diagonally we need to communicate over two interfaces. The second step of Section 3.2 may clarify this information transfer as already two iterations are necessary in two dimensions. Algorithm 1 shows the construction of tensor-product masks for a $d$-dimensional mesh.

**Remark 2.2.** The hierarchical nature of the integrated Legendre polynomials allows us to use inhomogeneous polynomial degrees between elements. If two connected elements have varying degrees in the tangential direction(s), we cannot trivially combine non-hierarchical shape functions to obtain globally continuous basis functions.

**Remark 2.3.** We can extend these ideas to unstructured meshes by considering that $d-1$ iterations may not be enough to reach all necessary elements and instead we could iterate until there are no more changes. Moreover, the operations on the tensor-product masks must consider that adjacent coordinate systems may not be aligned, which requires flipping some odd-degree polynomial shape functions.

**Algorithm 1** Construct compatible tensor-product masks for single $p$-finite element mesh.

```
// Compare the slices (n,1) and (n,0) of the nd-arrays A_0 and A_1 using operation op
// and overwrite them with the result. Data types are either bool or int.
void operateOnInterface(NdArray<Type>& A_0, NdArray<Type>& A_1, int n, Op op)
{
    // Extract sizes of arrays in each direction without entry at index n
    Vector s_0 = sizes(A_0);
    Vector s_1 = sizes(A_1);

    Vector s_0^n = removeEntry(s_0, n);
    Vector s_1^n = removeEntry(s_1, n);

    // Compute element-wise maximum s^n of s_0^n and s_1^n and iterate over all
    // entries in the Cartesian product [0, s^n(0)] × ... × [0, s^n(d - 2)]
    for(Vector i^n in productIndices(maxarray(s_0^n, s_1^n)))
    {
        Vector i_0 = insertEntry(i^n, n, 1); // second face of first element
        Vector i_1 = insertEntry(i^n, n, 0); // first face of second element

        bool e_0 = indexExists(i_0, s_0);
        bool e_1 = indexExists(i_1, s_1);

        Type v_0 = A_0(i_0) if e_0 else noValue(Type);
        Type v_1 = A_1(i_1) if e_1 else noValue(Type);

        Type r = op(v_0, v_1);

        if(v_0 ≠ r and not e_0) resize(A_0, i_0);
        if(v_1 ≠ r and not e_1) resize(A_1, i_1);

        if(v_0 ≠ r) A_0(i_0) = r;
        if(v_1 ≠ r) A_1(i_1) = r;
    }
}

// Compare array slices of all cells with second neighbor if it exists
void operateOnInterfaces(NdArrayList<Type>& A, Neighbors N, Op op)
{
    for(int a from 0 to d - 1)          // Loop over coordinate axes
        for(int i from 0 to size(N))    // Loop over cell indices
            if(int N_a = N(i,a,1); N_a ≠ −1)
                operateOnInterface(A(i), A(N_a), a, op);
}

MaskList createPfemMasks(Neighbors N, Degrees p, String strategy)
{
    MaskList M;

    // Activate full tensor-product (could use trunk space masks instead)
    for(Vector p_e in p)
        M.append(Mask(p_e + 1, true));

    Op op = logicalAnd if strategy == "min_degree" else logicalOr;

    // Recover compatibility
    for(int it from 0 to d - 2)
        operateOnInterfaces(M, N, op);

    return M;
}
```

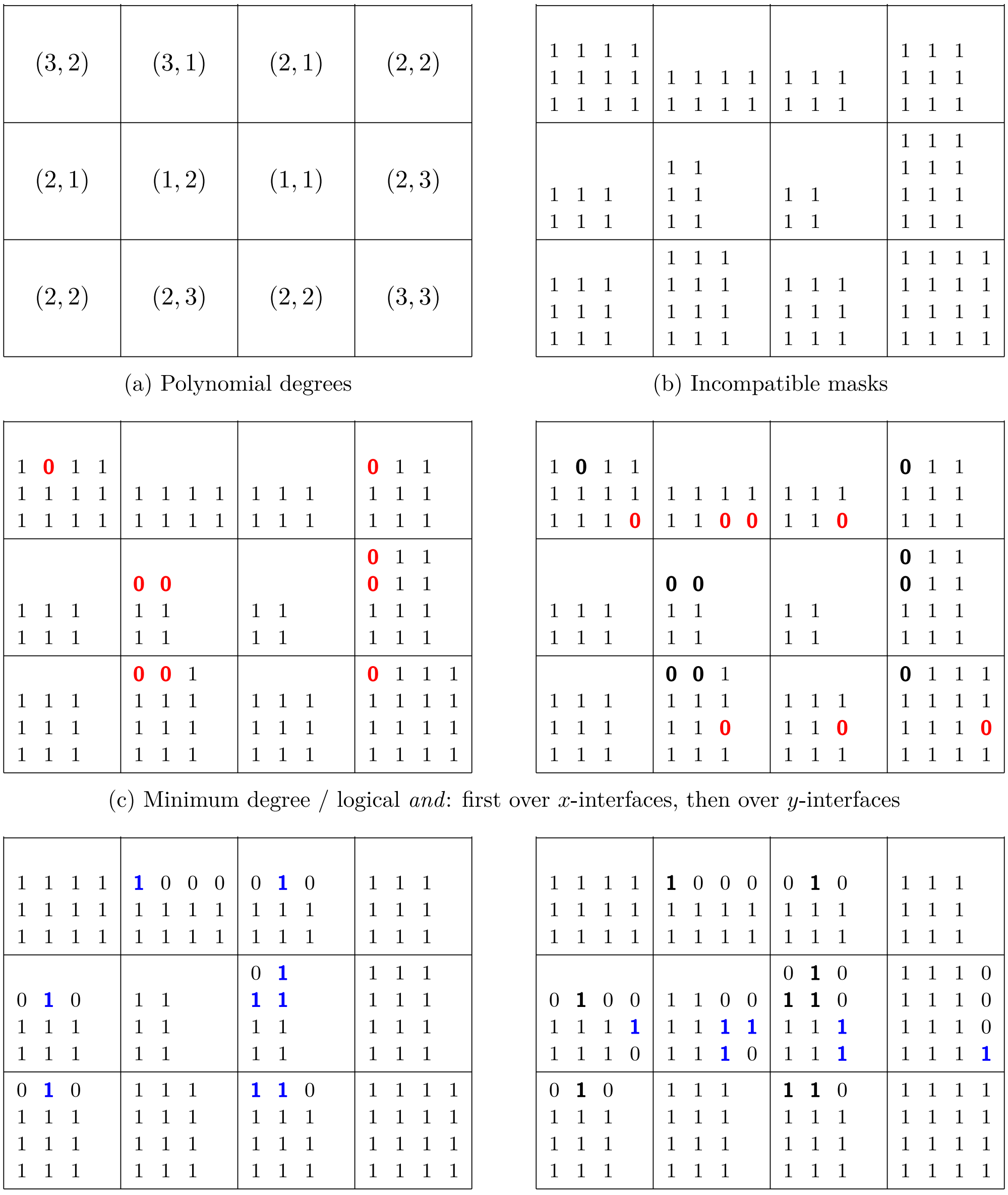

(a) Polynomial degrees

(b) Incompatible masks

(c) Minimum degree / logical *and*: first over $x$-interfaces, then over $y$-interfaces

(d) Maximum degree / logical *or*: first over $x$-interfaces, then over $y$-interfaces

**Fig. 5.** Construct compatible tensor-product masks for a $p$-finite element mesh.

## 2.3. Location matrices

The tensor-product masks we constructed in the previous section determine the active shape functions for all elements, such that they have matching counterparts across all interfaces. The remaining task is to connect sets of shape functions from different elements into basis functions by assigning them to the same global index. To store this data, we introduce a *location matrix* for each tensor-product mask with the same size that contains the global basis function index for each active entry and some value, e.g., $-1$ or the maximum representable number for inactive entries. In higher dimensions, the location matrices are $d$-dimensional arrays containing the basis function indices.

**Algorithm 2** Construct location matrices for a single $p$-finite element mesh.

```
[LocationMatrixList, int] initializeGlobalIndices(MaskList M)
{
    LocationMatrixList G; int n_ids = 0;

    for(int i from 0 to size(M) - 1)
    {
        Mask M_i = M(i);
        LocationMatrix G_i(sizes(M_i));

        // Loop over all tensor-product entries (d nested loops)
        for(Vector j in productIndices(sizes(M_i)))
            G_i(j) = n_ids++ if M_i(j) else -1;

        G.append(G_i);
    }

    return [G, n];
}

void removeUnassignedIndices(LocationMatrixList& G, int n_ids)
{
    Vector exists(n_ids, 0);
    Vector map(n_ids, -1);

    // Activate entries in exists vector
    for(LocationMatrix G_i in G)
        for(Vector j in productIndices(sizes(G_i)))
            if(G_i(j) ≠ -1)
                exists(G_i(j)) = 1;

    int n_new = 0;

    // Count active dofs upwards into map
    for(int i from 0 to nids - 1)
        map(i) = n_new++ if exists(i) else −1;

    // Reassign indices in location matrices
    for(LocationMatrix& G_i in G)
        for(Vector j in productIndices(sizes(G_i)))
            if(G_i(j) ≠ -1)
                G_i(j) = map(G_i(j))
}

// Operation for connecting shapes that returns first id
int returnFirstIndex(int i_0, int i_1)
{
    return i_0;
}

// Create location matrix with global ids for each tensor-product mask
void createPfemLocationMatrices(MaskList M, Neighbors N)
{
    LocationMatrixList G, int n_ids = initializeGlobalIndices(M);

    operateOnInterfaces(G, N, returnFirstIndex);
    removeUnassignedIndices(G, n_ids);

    return G;
}
```

| | | | |
|---|---|---|---|
| 17 -1 22 24<br>16 19 21 23<br>15 18 20 -1 | 39 41 42 43<br>38 40 -1 -1 | 57 59 60<br>56 58 -1 | -1 89 92<br>86 88 91<br>85 87 90 |
| 10 12 14<br>9 11 13 | -1 -1<br>35 37<br>34 36 | 53 55<br>52 54 | -1 80 84<br>-1 79 83<br>76 78 82<br>75 77 81 |
| 2 5 8<br>1 4 7<br>0 3 6 | -1 -1 33<br>27 30 32<br>26 29 -1<br>25 28 31 | 46 49 51<br>45 48 -1<br>44 47 50 | -1 67 71 74<br>63 66 70 73<br>62 65 69 -1<br>61 64 68 72 |

(a) Assign unique indices

| | | | |
|---|---|---|---|
| 17 -1 22 24<br>16 19 21 23<br>15 18 20 -1 | 19 41 42 43<br>18 40 -1 -1 | 41 59 60<br>40 58 -1 | -1 89 92<br>59 88 91<br>58 87 90 |
| 10 12 14<br>9 11 13 | -1 -1<br>12 37<br>11 36 | 37 55<br>36 54 | -1 80 84<br>-1 79 83<br>55 78 82<br>54 77 81 |
| 2 5 8<br>1 4 7<br>0 3 6 | -1 -1 33<br>5 30 32<br>4 29 -1<br>3 28 31 | 30 49 51<br>29 48 -1<br>28 47 50 | -1 67 71 74<br>49 66 70 73<br>48 65 69 -1<br>47 64 68 72 |

(b) Connect in $x$

| | | | |
|---|---|---|---|
| 17 -1 22 24<br>16 19 21 23<br>10 12 14 -1 | 19 41 42 43<br>12 37 -1 -1 | 41 59 60<br>37 55 -1 | -1 89 92<br>59 88 91<br>55 78 82 |
| 10 12 14<br>1 4 7 | -1 -1<br>12 37<br>4 29 | 37 55<br>29 48 | -1 80 84<br>-1 79 83<br>55 78 82<br>48 65 69 |
| 2 5 8<br>1 4 7<br>0 3 6 | -1 -1 33<br>5 30 32<br>4 29 -1<br>3 28 31 | 30 49 51<br>29 48 -1<br>28 47 50 | -1 67 71 74<br>49 66 70 73<br>48 65 69 -1<br>47 64 68 72 |

(c) Connect in $y$

| | | | |
|---|---|---|---|
| 13 -1 16 18<br>12 14 15 17<br>9 10 11 -1 | 14 26 27 28<br>10 25 -1 -1 | 26 35 36<br>25 34 -1 | -1 55 57<br>35 54 56<br>34 48 51 |
| 9 10 11<br>1 4 7 | -1 -1<br>10 25<br>4 20 | 25 34<br>20 30 | -1 50 53<br>-1 49 52<br>34 48 51<br>30 38 42 |
| 2 5 8<br>1 4 7<br>0 3 6 | -1 -1 24<br>5 21 23<br>4 20 -1<br>3 19 22 | 21 31 33<br>20 30 -1<br>19 29 32 | -1 40 44 47<br>31 39 43 46<br>30 38 42 -1<br>29 37 41 45 |

(e) Eliminate unassigned indices using (d)

| | | | | | | | | | | | | | | | | | | | | | | | | |
|---|---|---|---|---|---|---|---|---|---|---|---|---|---|---|---|---|---|---|---|---|---|---|---|---|
| old index: | 0 | 1 | 2 | 3 | 4 | 5 | 6 | 7 | 8 | 9 | 10 | 11 | 12 | 13 | 14 | 15 | 16 | 17 | 18 | 19 | 20 | 21 | 22 | 23 |
| exists: | 1 | 1 | 1 | 1 | 1 | 1 | 1 | 1 | 1 | 0 | 1 | 0 | 1 | 0 | 1 | 0 | 1 | 1 | 0 | 1 | 0 | 1 | 1 | 1 |
| new index: | 0 | 1 | 2 | 3 | 4 | 5 | 6 | 7 | 8 | -1 | 9 | -1 | 10 | -1 | 11 | -1 | 12 | 13 | -1 | 14 | -1 | 15 | 16 | 17 |
| old index: | 24 | 25 | 26 | 27 | 28 | 29 | 30 | 31 | 32 | 33 | 34 | 35 | 36 | 37 | 38 | 39 | 40 | 41 | 42 | 43 | 44 | 45 | 46 | 47 |
| exists: | 1 | 0 | 0 | 0 | 1 | 1 | 1 | 1 | 1 | 1 | 0 | 0 | 0 | 1 | 0 | 0 | 0 | 1 | 1 | 1 | 0 | 0 | 0 | 1 |
| new index: | 18 | -1 | -1 | -1 | 19 | 20 | 21 | 22 | 23 | 24 | -1 | -1 | -1 | 25 | -1 | -1 | -1 | 26 | 27 | 28 | -1 | -1 | -1 | 29 |
| old index: | 48 | 49 | 50 | 51 | 52 | 53 | 54 | 55 | 56 | 57 | 58 | 59 | 60 | 61 | 62 | 63 | 64 | 65 | 66 | 67 | 68 | 69 | 70 | 71 |
| exists: | 1 | 1 | 1 | 1 | 0 | 0 | 0 | 1 | 0 | 0 | 0 | 1 | 1 | 0 | 0 | 0 | 1 | 1 | 1 | 1 | 1 | 1 | 1 | 1 |
| new index: | 30 | 31 | 32 | 33 | -1 | -1 | -1 | 34 | -1 | -1 | -1 | 35 | 36 | -1 | -1 | -1 | 37 | 38 | 39 | 40 | 41 | 42 | 43 | 44 |
| old index: | 72 | 73 | 74 | 75 | 76 | 77 | 78 | 79 | 80 | 81 | 82 | 83 | 84 | 85 | 86 | 87 | 88 | 89 | 90 | 91 | 92 | | | |
| exists: | 1 | 1 | 1 | 0 | 0 | 0 | 1 | 1 | 1 | 0 | 1 | 1 | 1 | 0 | 0 | 0 | 1 | 1 | 0 | 1 | 1 | | | |
| new index: | 45 | 46 | 47 | -1 | -1 | -1 | 48 | 49 | 50 | -1 | 51 | 52 | 53 | -1 | -1 | -1 | 54 | 55 | -1 | 56 | 57 | | | |

(d) Index map

**Fig. 6.** Glue shape functions from Fig. 5(c) together into compatible basis functions.

We begin by assigning all shape functions a unique global index, as Fig. 6(a) shows. Then, we connect the shape functions across all element interfaces by copying the indices from one side to the other. Figs. 6(b) and 6(c) demonstrate this process by connecting the indices in the $x$- and then in the $y$-direction. This step leaves some indices without contributions from any elements, as their original shape functions have been assigned to a different index. We eliminate these "gaps" in the global numbering by constructing a new set of compressed indices in Fig. 6(e) that we obtain by applying the index map shown in Fig. 6(d). Algorithm 2 summarizes the creation of location matrices that finalizes the construction of our $p$-finite element basis.

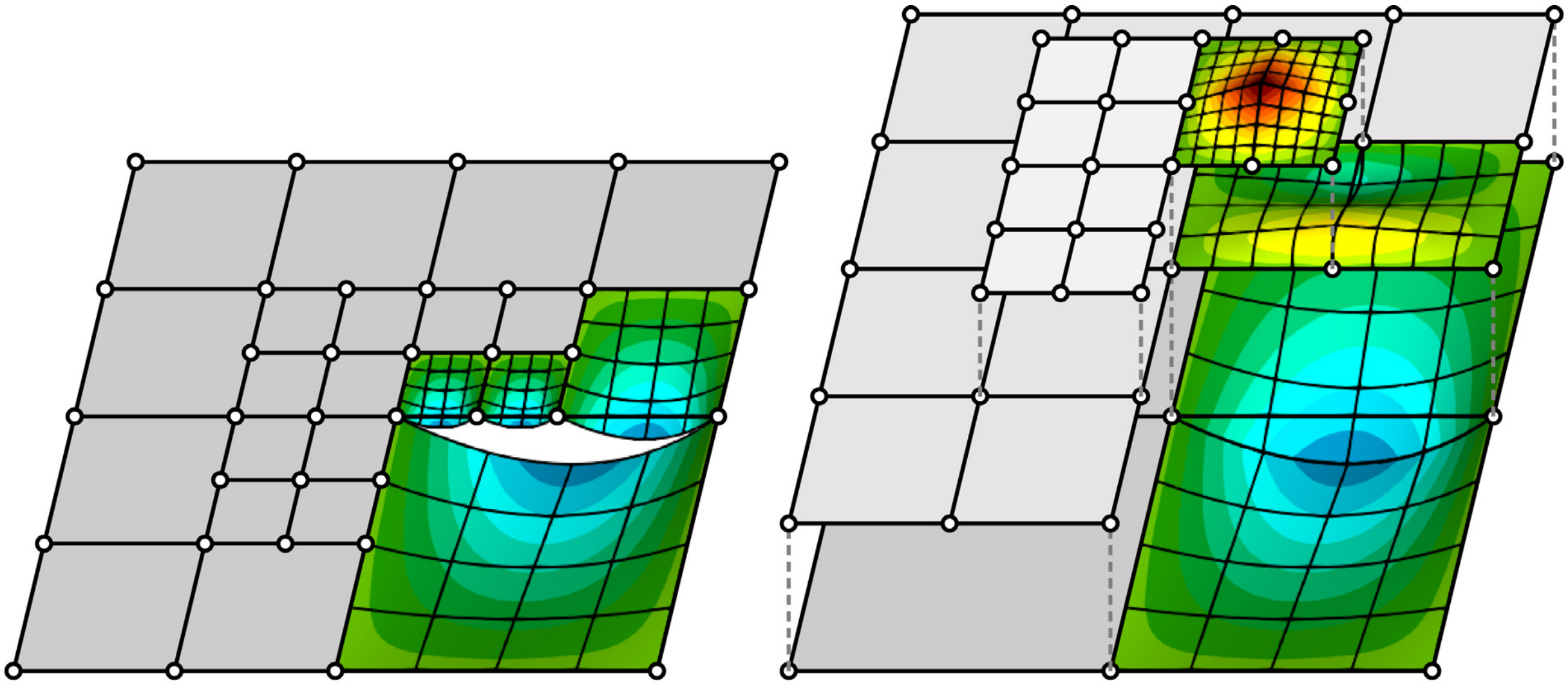

(a) Refine by replacement (b) Refine by superposition

**Fig. 7.** Comparison of $hp$-refinement strategies.
*Source:* Reprinted from [12].

**Remark 2.4.** By linearizing a location matrix and filtering the inactive entries ($-1$ values), we obtain a global index vector that is sometimes also called location map or element freedom table.

## 3. Multi-level $hp$-finite element basis

Refining the finite element discretization by only increasing the polynomial degree is often not the best strategy as we also need to reduce $h$ to recover exponential convergence for solutions with reduced regularity. Even for smooth solutions, it is often more efficient to start with several levels of $h$-refinement until the solution becomes smooth enough on each element, as $p$-refinement is more expensive per degree of freedom than $h$-refinement. However, replacing coarse finite elements by a set of finer ones introduces hanging nodes—element interfaces that are subdivided on one side, but not on the other. The presence of hanging nodes greatly complicates the construction of continuous basis functions as the shape functions of finer elements cannot be directly combined with shape functions on the coarser element. Instead, combinations of shape functions must be constrained and connected. The complexity of this process increases drastically with more spatial dimensions and higher degrees of mesh irregularity [6,7,9].

The multi-level $hp$-method circumvents the challenge of hanging nodes by keeping the coarse elements during the refinement. By allowing these coarse elements to support shape functions, we resolve cases where neighboring elements have a different refinement level. Instead of trying to connect the coarse shape functions to the finer ones, we simply activate the functions on the original coarse neighbor (that is now overlayed with finer elements). Using this idea, we combine shape functions from neighboring elements, essentially reducing the problem of constructing an $hp$-basis to creating multiple levels of $p$-finite element meshes. Fig. 7 compares the replacement and overlay strategies and shows how coarse shape functions are completed on the coarse parent of the finer neighbors.

The multi-level $hp$-method introduces two main challenges: we need to preserve linear independence and our overlay functions must be zero on internal boundaries. In the following, we address these challenges and show how we automatically resolve them with little implementation effort.

### 3.1. Hierarchical mesh data structure

We begin our refinement process with a base mesh, where we overlay cells marked for refinement by $2^d$ sub-cells from bisecting the original cell in each coordinate direction. Performing this process recursively on overlay cells

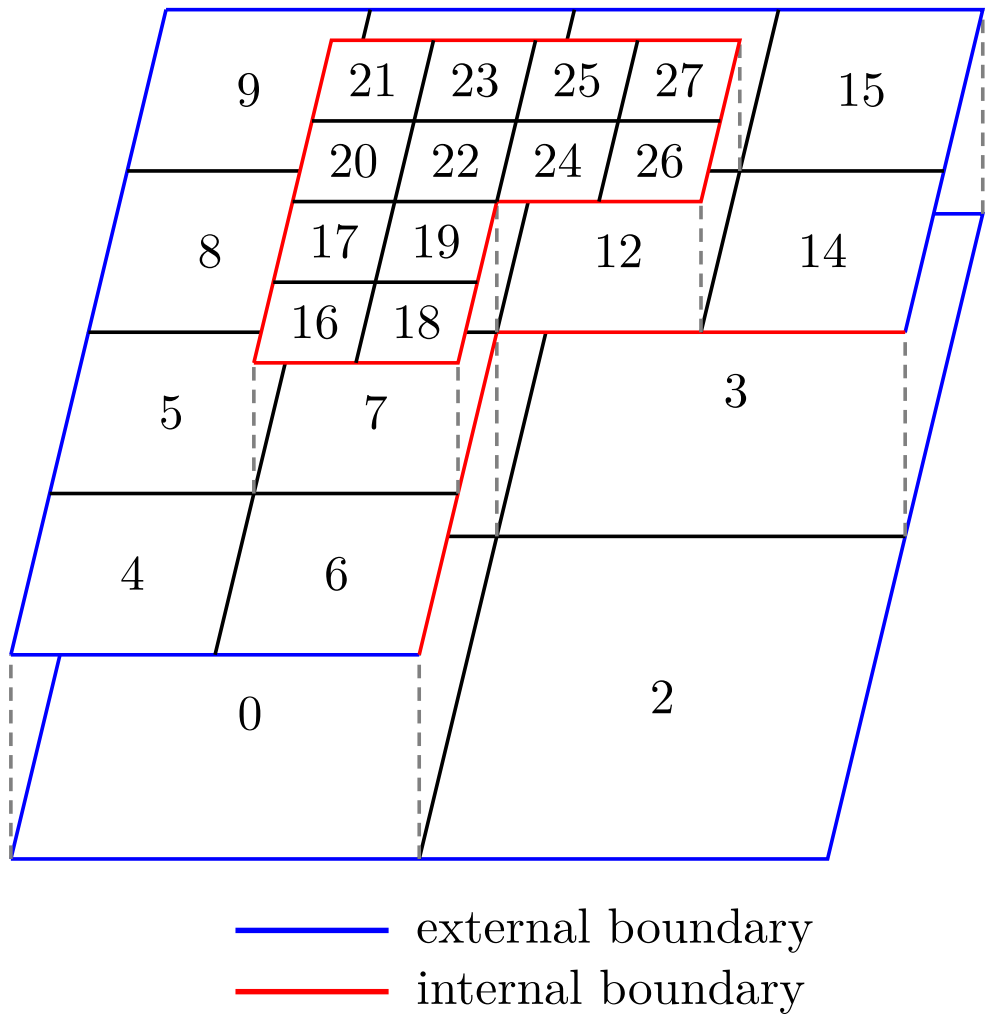

(a) Full cell numbering

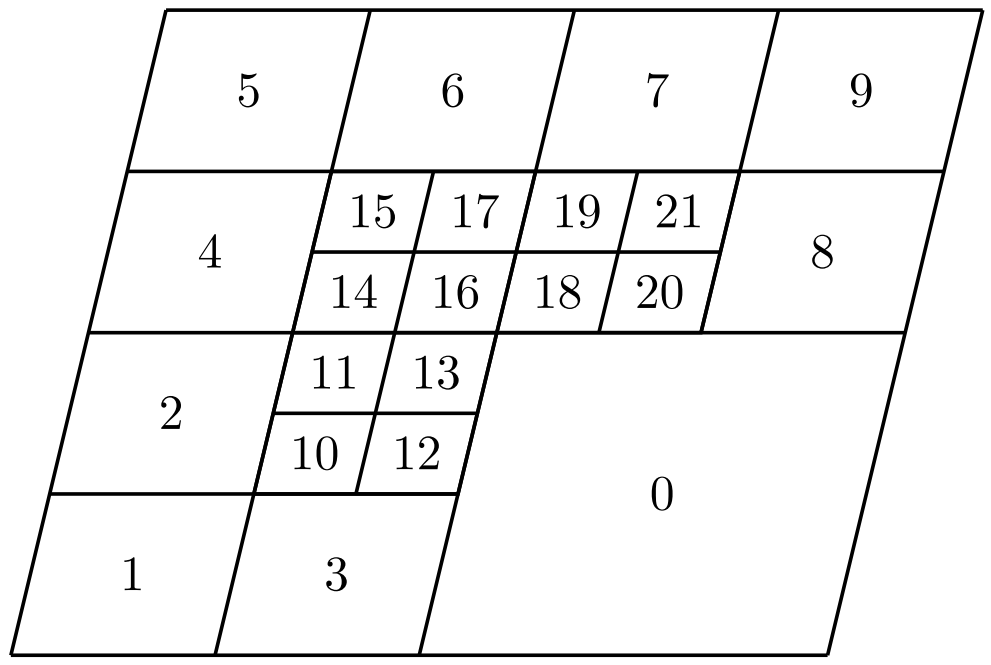

(b) Leaf cell (element) numbering

| Cell | $n_l$ | $n_r$ | $n_b$ | $n_t$ | $P$ | $R$ | $L$ |
|---|---|---|---|---|---|---|---|
| 0 | -1 | 2 | -1 | 1 | -1 | 0 | 0 |
| 1 | -1 | 3 | 0 | -1 | -1 | 0 | 0 |
| 2 | 0 | -1 | -1 | 3 | -1 | 0 | 1 |
| 3 | 1 | -1 | 2 | -1 | -1 | 0 | 0 |
| 4 | -1 | 6 | -1 | 5 | 0 | 1 | 1 |
| 5 | -1 | 7 | 4 | 8 | 0 | 1 | 1 |
| 6 | 4 | 2 | -1 | 7 | 0 | 1 | 1 |
| 7 | 5 | 2 | 6 | 10 | 0 | 1 | 0 |
| 8 | -1 | 10 | 5 | 9 | 1 | 1 | 1 |
| 9 | -1 | 11 | 8 | -1 | 1 | 1 | 1 |
| 10 | 8 | 12 | 7 | 11 | 1 | 1 | 0 |
| 11 | 9 | 13 | 10 | -1 | 1 | 1 | 1 |
| 12 | 10 | 14 | 2 | 13 | 3 | 1 | 0 |
| 13 | 11 | 15 | 12 | -1 | 3 | 1 | 1 |
| 14 | 12 | -1 | 2 | 15 | 3 | 1 | 1 |
| 15 | 13 | -1 | 14 | -1 | 3 | 1 | 1 |
| 16 | 5 | 18 | 6 | 17 | 7 | 2 | 1 |
| 17 | 5 | 19 | 16 | 20 | 7 | 2 | 1 |
| 18 | 16 | 2 | 6 | 19 | 7 | 2 | 1 |
| 19 | 17 | 2 | 18 | 22 | 7 | 2 | 1 |
| 20 | 8 | 22 | 17 | 21 | 10 | 2 | 1 |
| 21 | 8 | 23 | 20 | 11 | 10 | 2 | 1 |
| 22 | 20 | 24 | 19 | 23 | 10 | 2 | 1 |
| 23 | 21 | 25 | 22 | 11 | 10 | 2 | 1 |
| 24 | 22 | 26 | 2 | 25 | 12 | 2 | 1 |
| 25 | 23 | 27 | 24 | 13 | 12 | 2 | 1 |
| 26 | 24 | 14 | 2 | 27 | 12 | 2 | 1 |
| 27 | 25 | 14 | 26 | 13 | 12 | 2 | 1 |

(c) Data structure

**Fig. 8.** Essential topological data: left, right, bottom and top neighbors $n_l$, $n_r$, $n_b$, $n_t$, parent index $P$, refinement level $R$, and a flag $L$ marking the leaf cells. (For interpretation of the references to color in this figure legend, the reader is referred to the web version of this article.)

yields a $2^d$-ary refinement tree. The choice to refine can use an error indicator or a priori knowledge of the solution. All cells within our refinement tree can support basis functions, although we choose finite elements only as leaf cells.

A central aspect of our framework is to identify a lightweight data structure that can store all essential information for constructing an *hp*-basis. The algorithms we introduce in this section need for each cell: the parent, the refinement level, the neighbors on the same or a lower refinement level, and the information whether the cell is a leaf. We store the neighbors in an $(n_{\text{cell}}, d, 2)$ sized array (denoted as $N$) and use simple one-dimensional arrays with one entry per cell for the other relations. By using the parent cell index $-1$ for cells without parent we identify the roots of the refinement tree. In our neighbor array, we distinguish between *internal interfaces* (neighbor exists and has same refinement level), *internal boundaries* (neighbor exists and has lower refinement level) and *external boundaries* (neighbor does not exist, face is part of the domain boundary). We again assign the *no cell* value for external boundaries and the index of the coarser leaf cell for internal boundaries. Fig. 8 sketches a simple example mesh with two levels of refinement.

**Remark 3.1.** Our internal boundaries are the multi-level *hp* equivalent of the hanging nodes from classical *hp*-finite elements. The treatment of hanging nodes introduces constraints to the shape functions on both sides, however, in the multi-level $hp$-method, we simply impose homogeneous Dirichlet conditions on internal boundaries (by removing active shape functions on them).

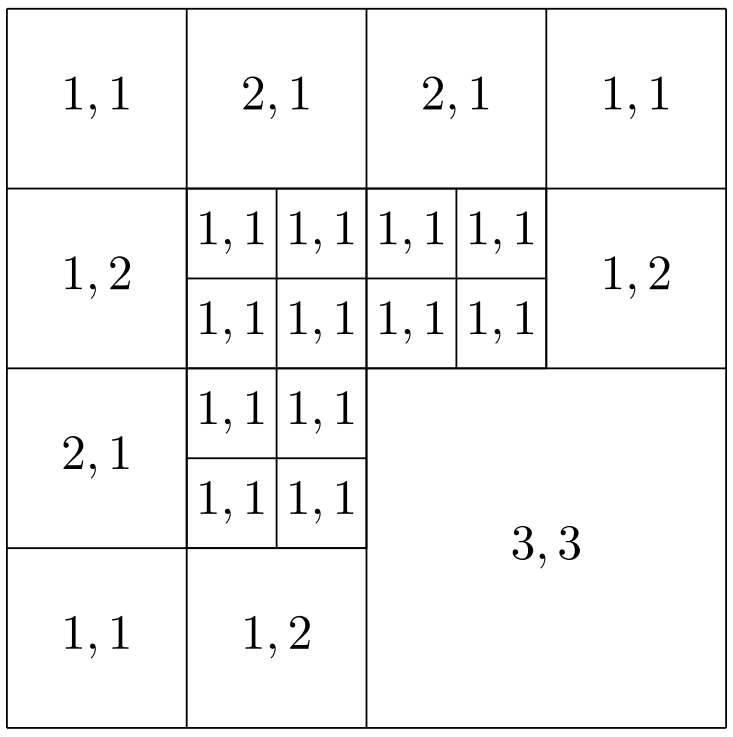

**Fig. 9.** Assign polynomial degrees to each finite element (Fig. 8(b)).

**Remark 3.2.** One can exploit the redundancy in our data structure to reduce the amount of storage that is needed for the meshes. For example, the refinement levels and the leaf flags can be trivially constructed on-the-fly from the parent indices. However, additional information like the geometric mapping function of each cell may be needed on top.

**Remark 3.3.** Whether $N$ denotes the shape functions or the neighbor array depends on the context.

### *3.2. Tensor-product masks*

We now construct tensor-product masks for all cells (leaves and non-leaves) in the refinement tree introduced in the previous section. These tensor-product masks activate the right shape functions on each cell, such that they can combine into globally continuous, linear independent basis functions that are complete up to the selected polynomial degrees (at least in the element interior). These polynomial degrees are only assigned to the leaves (finite elements), and they can again vary in the local directions. For our example mesh, we choose the polynomial degrees shown in Fig. 9. The following algorithm extends the ideas of Section 2.2 to construct the multi-level $hp$ tensor-product masks in four steps (Figs. 10–13 and Algorithm 3).

**1. Initialize leaf masks according to their polynomial degree and non-leaf masks empty**. By activating the complete tensor-product or *trunk space* (see the Appendix) only on the leaves, we naturally preserve linear independence as currently none of the parents contribute any shape functions. However, the resulting masks in Fig. 10 cannot be combined into continuous basis functions yet as they are not compatible along the cell interfaces. Moreover, the hierarchical nature of the integrated Legendre shape functions allows to remove certain modes from the tensor-product that do not contribute to the order of convergence. The result is commonly referred to as the trunk space. The Appendix derives the corresponding initial tensor-product masks that can be applied in this step instead.

**2. Restore interface compatibility using logical *or* operations in $d$ iterations.** We loop over the internal interfaces and activate entries in the tensor-product masks if the corresponding entry in the neighbor is active. Internal interfaces bound two cells of the same refinement level, which makes this step equivalent to treating the cells of each refinement level as separate $p$-finite element meshes and applying the algorithm of Section 2.2 to them independently. This serves two purposes: first, we resolve incompatible polynomial degrees by the *maximum degree strategy*, and second, we activate the interface shape functions on the parent elements of finer neighbors. Fig. 11 shows how the shape functions on the bottom right cell of level 0 can now be completed on the neighboring cells of the same level, although these cells have been refined. Two iterations of our algorithm are required to activate all four linear shape functions that together form the central hat basis function on level 0.

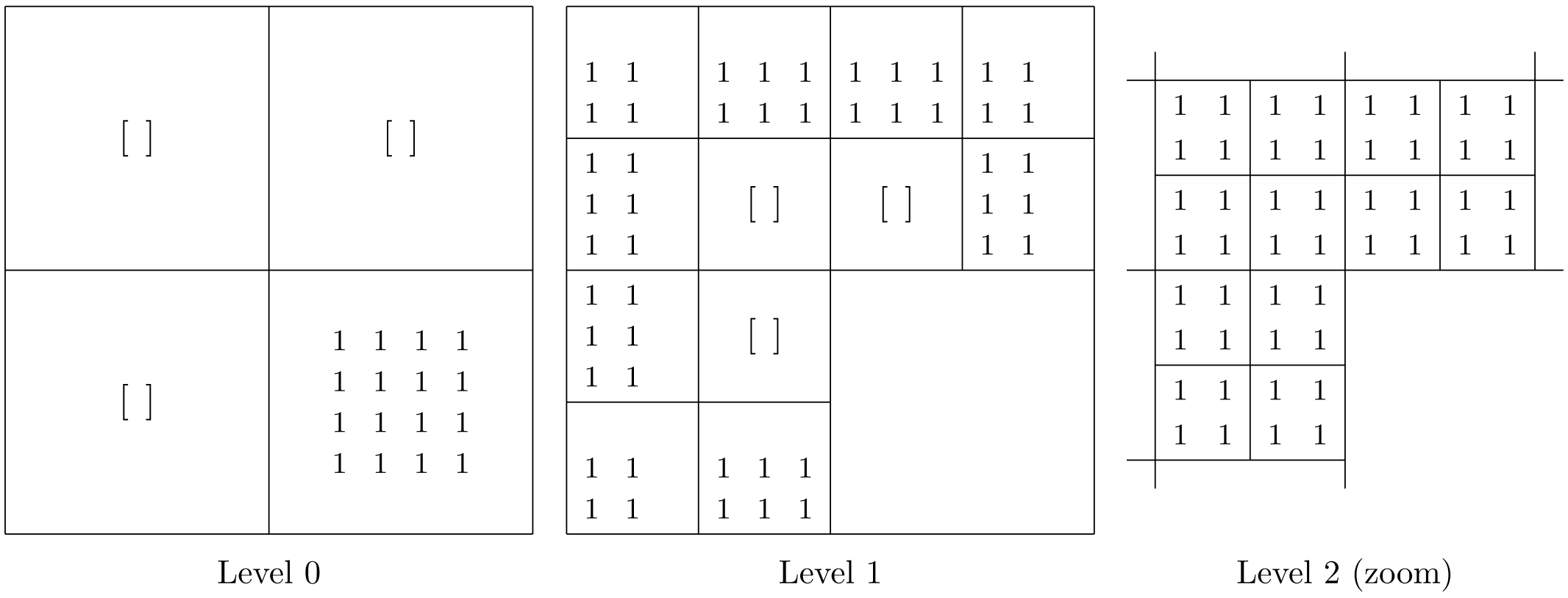

**Fig. 10.** Initial tensor-product masks corresponding to degrees in Fig. 9.

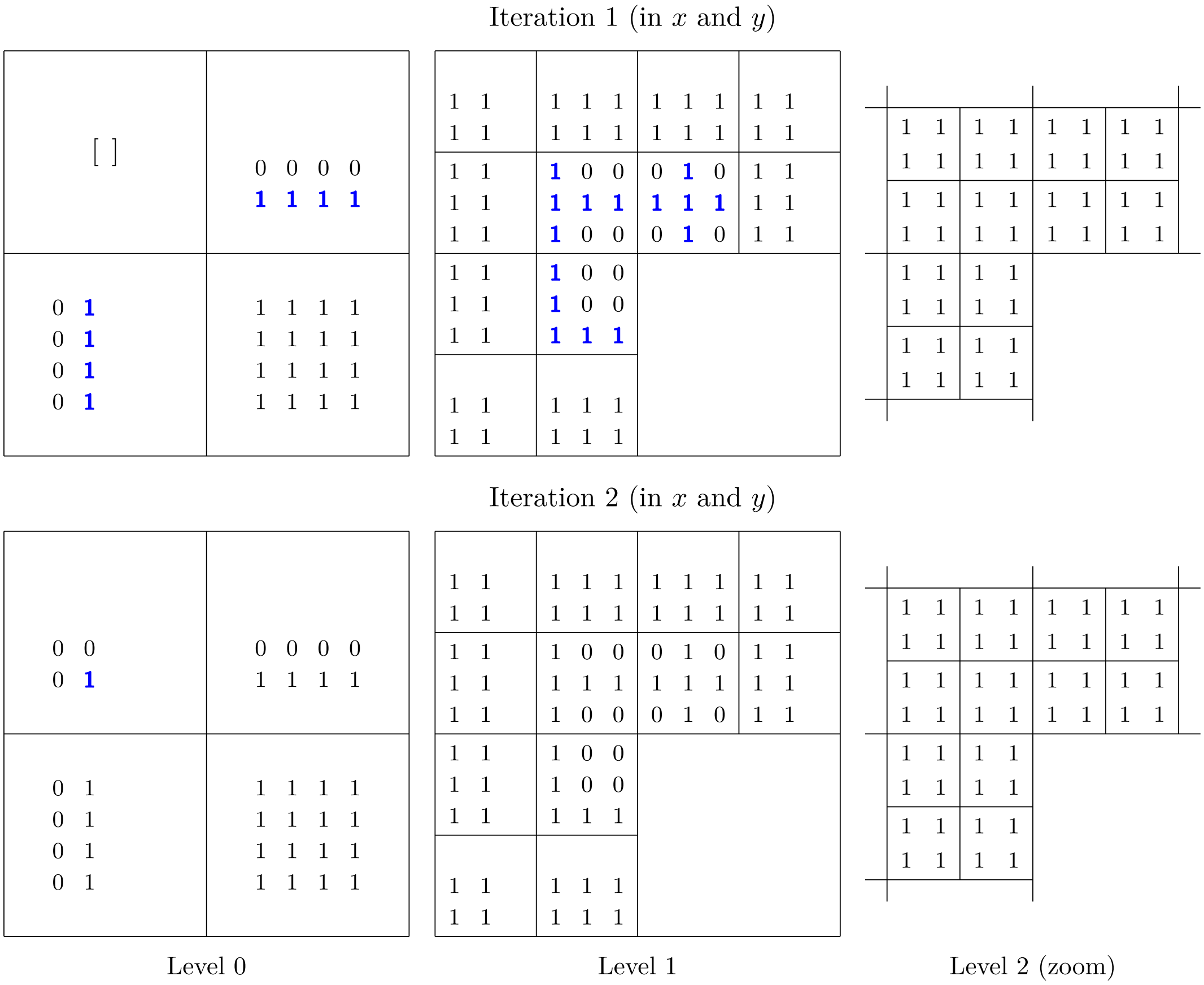

**Fig. 11.** Restore continuity after initial activation.

**3. Deactivate shape functions on internal boundaries.** We finalize the construction of our multi-level $hp$-basis by ensuring that our basis functions are zero on the internal boundaries of our overlays (i.e., imposing homogeneous Dirichlet conditions there). Thus, we deactivate the entries in the tensor-product masks for shape functions that are non-zero on faces with a coarser neighbor. As Fig. 12 shows, we remove the shape functions on levels 1 and 2 that are active on the internal boundaries to the bottom right cell of level 0. Note, that the shape functions we deactivate in this step (on levels 1 and 2) are active on one of the parent cells (level 0) from the previous step. Moreover, we do not modify external boundaries as there are no hanging nodes and we want the finest resolution there (as already initialized in step 1).

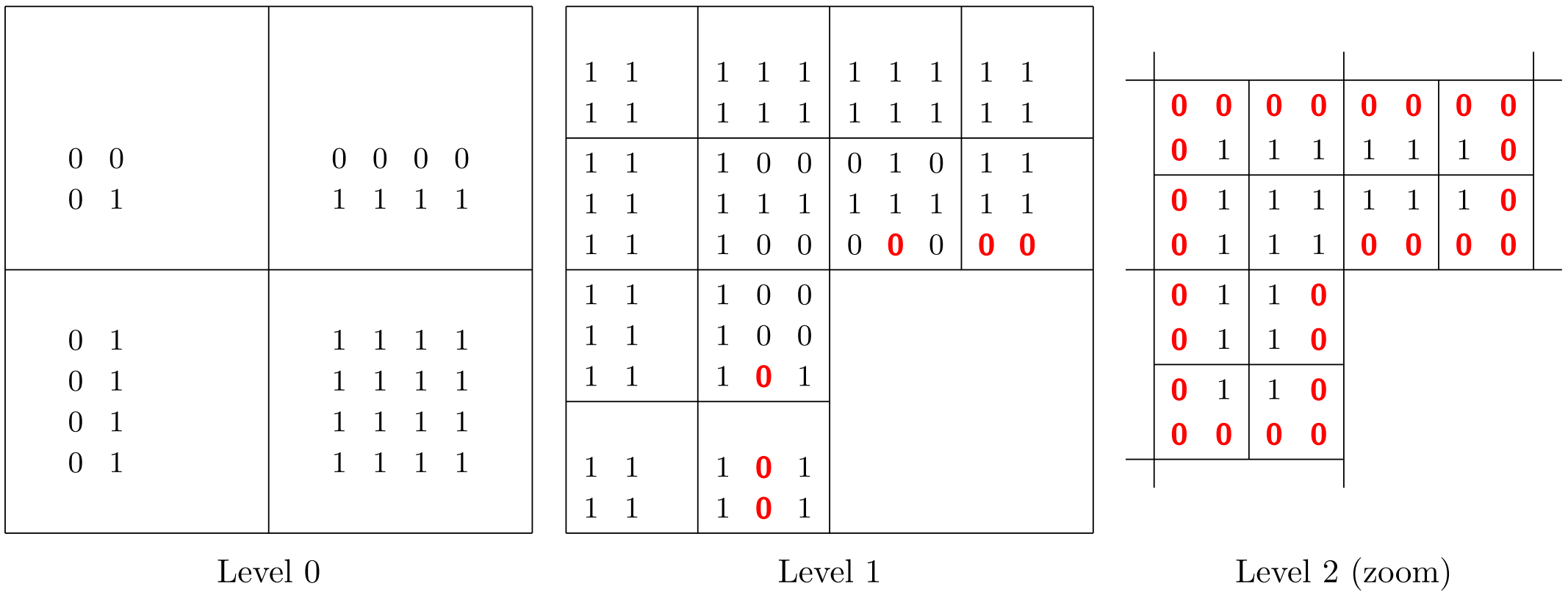

**Fig. 12.** Deactivate internal boundaries (step 3).

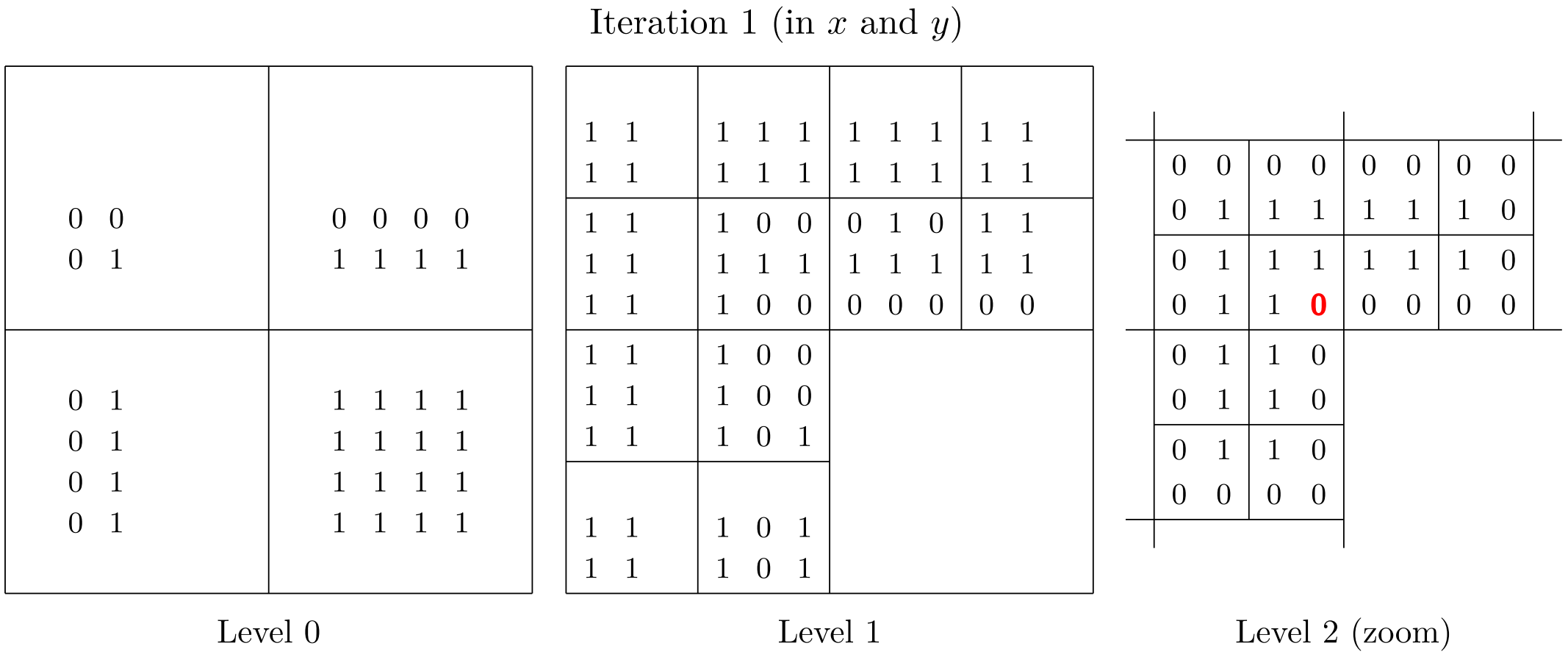

**Fig. 13.** Communicate zero Dirichlet conditions to neighbors (step 4).

**4. Restore interface compatibility using logical *and* operations in $d - 1$ iterations.** The basis functions on internal boundaries often receive contributions from additional elements. For example, on level 2 of Fig. 12, the basis function associated to the node in the central corner is built from three linear shape functions supported on the three surrounding elements (cells 19, 22, and 24), but only two of them are deactivated (on cells 19 and 24, as only they are directly connected to the internal boundary). We deactivate all contributions by repeating step 2 of our algorithm, but with one less repetition and with logical *and* operations to deactivate interface shape functions if they are not active in the neighbor, which only affects the changes from step 3, as all other interfaces are compatible after step 2. Fig. 13 shows the final set of tensor-product masks for our example mesh.

At this point, we can evaluate the basis functions in each finite element by first evaluating the integrated Legendre polynomials in each direction and then computing the active entries within the tensor-product according to the corresponding mask. If needed, we map the derivatives into the global coordinate system (physical space). We start on a leaf cell and repeat this process for each parent cell in the hierarchy towards the root cell. With each new level,

**Algorithm 3** Construct multi-level *hp* tensor-product masks and location matrices.

```
// Replaces p-fem version; requires neighbor to exist with same level
void operateOnInterfaces(NdArrayList<Type>& A, Neighbors N, Levels L, Op op)
{
    for(int a from 0 to d - 1)              // For d coordinate axes
        for(int i from 0 to size(A) - 1)    // For all cells
            if(int N_a = N(i, a, 1); N_a ≠ −1 and L(i) == L(N_a))
                operateOnInterface(A(i), A(N_a), a, op);
}

// Deactivate first (f == 0) or second (f == 1) slice normal to axis n
void deactivateSlice(Mask& M_i, int n, int f)
{
    Vector s^n = removeEntry(sizes(M_i), n);

    for(Vector i^n in productIndices(s^n))
        M_i(insertEntry(i^n, n, f)) = false;
}

// Deactivates boundary functions if neighbor exists with different level
void deactivateOnInternalBoundaries(MaskList& M, Neighbors N, Levels L)
{
    for(int a from 0 to d - 1)             // For d coordinate axes
        for(int f from 0 to 1)             // For two faces
            for(int i from 0 to size(M))   // For all cells
                if(int N_a = N(i, a, 1); N_a ≠ −1 and L(i) != L(N_a))
                    deactivateSlice(M(i), a, f);
}

MaskList createMlhpMasks(Neighbors N, Levels L, Isleaf E, Degrees p)
{
    MaskList M;

    // Step 1: Initialize leaf masks (could use trunk space here instead)
    l = 0
    for(bool i from 0 to size(E) - 1)
        M.append(Mask(p(l++) + 1, true) if E(i) else Mask(0));

    // Step 2: Restore compatibility by activating on neighbors
    for(int it from 0 to d - 1)
        operateOnInterfaces(M, N, L, logicalOr);

    // Step 3: Impose Dirichlet on internal boundaries
    deactivateOnInternalBoundaries(M, N, L);

    // Step 4: Restore compatibility by deactivating on neighbors
    for(int it from 0 to d - 2)
        operateOnInterfaces(M, N, L, logicalAnd);

    return M;
}

void createMlhpLocationMatrices(MaskList M, Neighbors N, Levels L)
{
    LocationMatrixList G, int n_ids = initializeGlobalIndices(M);
    operateOnInterfaces(G, N, L, returnFirstIndex);
    removeUnassignedIndices(G, n_ids);

    return G;
}
```

**Algorithm 4** Evaluate multi-level *hp*-basis functions.

```
// Evaluate multi-level hp basis at local coordinates r of cell i, differentiated
// k(a) times in direction a. For example, k = (1, 0) evaluates ∂N_j/∂r_0 in 2D and
// k = (0, 1, 1) evaluates ∂^2 N_j/(∂y ∂r_2) in 3D. Note, that this function only maps
// the derivatives into the leaf coordinate system, not yet into global space.
Vector evaluate(Mesh mesh, MaskList M, int i_l, Coords r, Diff k)
{
    Vector N = [ ];   // shape functions
    double n_L = 0;   // number of levels

    // Loop from leaf to root: start with full index of leaf and update to
    // parent in each new iteration until current cell has no parent (-1)
    for(int i = i_l; i ≠ -1; i = mesh.parent(i))
    {
        Mask M_i = M(i);
        Vector s = sizes(M_i);

        if(product(s) ≠ 0)
        {
            // Evaluate integrated Legendre polynomials
            VectorList I;

            for(int a from 0 to d - 1)
                I.append(integratedLegendre(r(a), s(a) - 1, k(a)));

            // Append active entries in tensor-product to N
            for(Vector α in productIndices(s))
            {
                if(M_i(α))
                {
                    // Map derivatives to leaf
                    double N_α = power(1/2, n_L * sum(k));

                    // Multiply in directions j with polynomial at index α(j),
                    // differentiated k(j) times, evaluated at r(j)
                    for(int j from 0 to d - 1)
                        N_α = N_α * I(j)(α(j))

                    N.append(N_α);
                }
            }
        }

        // Map coordinate to parent
        for(int a from 0 to d - 1)
        {
            int N_a = mesh.neighbor(i, a, 0);

            double c = -1 if N_a == -1 or mesh.parent(N_a) != mesh.parent(i) else 1

            r(a) = (r(a) + c) / 2;
        }

        n_L = n_L + 1;
    }

    return N;
}
```

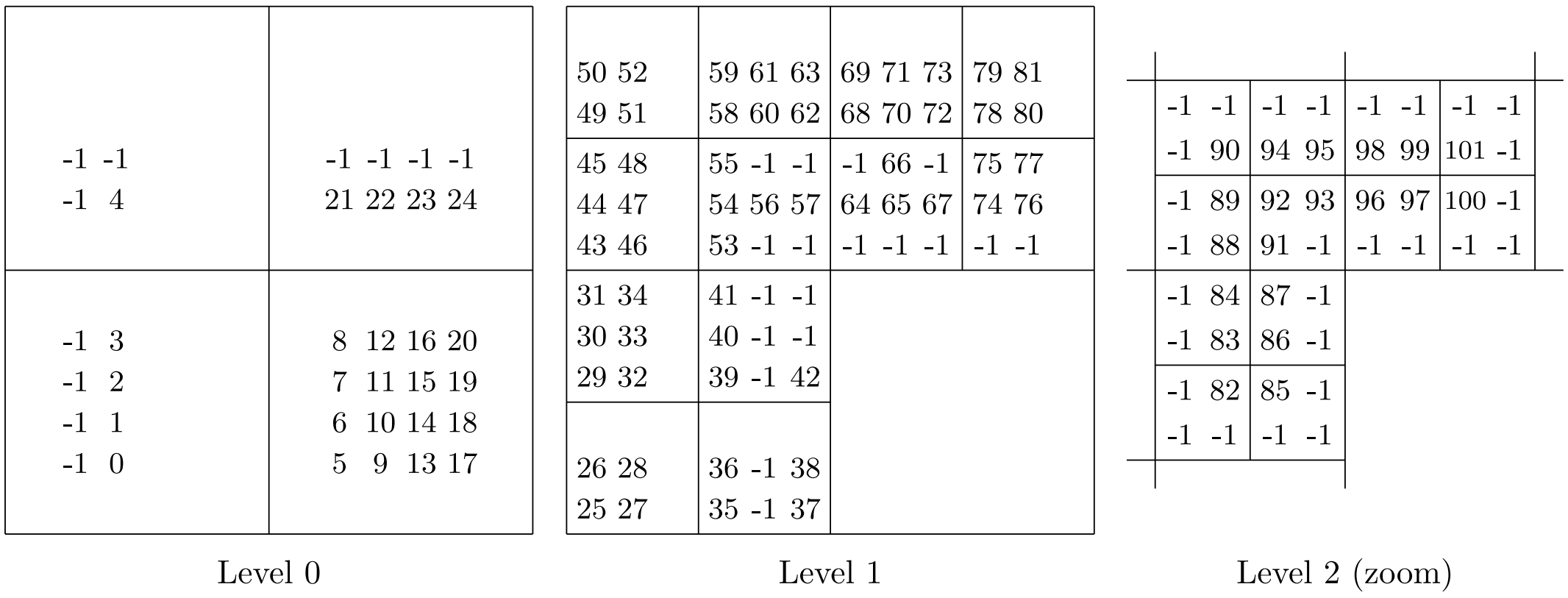

**Fig. 14.** Initialize global ids independently.

we map the evaluation coordinates accordingly. Algorithm 4 shows the basis function evaluation for a multi-level $hp$-element.

**Remark 3.4.** We introduce the tensor-product masks as individual tensors. Alternatively, one could first determine bounds for their size and then allocate them once potentially in a contiguous fashion.

### *3.3. Location matrices*

We connect the active shape functions on the cells of our refinement tree from the previous section to global multi-level $hp$-basis functions by using the ideas of Section 2.2. We treat each refinement level like a separate $p$-finite element mesh and construct the location matrices using the same steps as in Section 2.2. In other words: we connect the shape functions only across internal interface and not across interfaces to coarser neighbors (internal boundaries). Fig. 14 shows the initial location matrices corresponding to the tensor-product masks from Fig. 13, where we assign each shape function a unique index. These are then connected in both coordinate directions in Fig. 15 and the resulting gaps in the global numbering are eliminated in Fig. 16 by applying an index map like in Fig. 6(d). The last function of Algorithm 3 constructs the location matrices for a multi-level $hp$-basis.

**Remark 3.5.** The chosen example features only two levels of refinement. The algorithms, however, apply to arbitrary nested space-tree partitions. We only distinguish between leaf- and non-leaf cells and between different types of cell boundaries.

**Remark 3.6.** In practice, one might encapsulate the topological and geometric data from Section 3.1 into a mesh data type and encapsulate the results of Sections 3.2–3.3 into an $hp$-basis data type.

### *3.4. Simulation workflow*

In a computation with multi-level $hp$-finite elements, the algorithms from Sections 3.2 and 3.3 directly follow the mesh creation to prepare the tensor-product masks and location matrices. Both are defined on all cells of the refinement tree. In theory, it is possible to consider each cell a finite element and assemble their contributions separately. However, the volumetric coupling of the basis functions through the hierarchy renders this approach impractical as the interactions between all the cells in the tree must be integrated. Moreover, the concept of overlapping finite elements is not very transparent for someone concerned with practical applications who may not be familiar with the details of the method. Instead, we consider only the leaves as finite elements and when evaluating the shape functions we simply append the contributions of all parent elements by traversing the hierarchy and mapping the evaluation coordinates accordingly. Therefore, from ”outside”, the multi-level $hp$-basis consists of the non-overlapping leaves, where some of them evaluate more shape functions than others (for example in transition zones with many hanging nodes). Hiding the hierarchical nature of the basis (that internally still exists) allows us to handle our multi-level *hp*-basis like any other finite element method and use standard algorithms. We also concatenate the global basis function indices in all location matrices through the hierarchy to a single location map for each finite element. Fig. 17 shows these location maps for our example.

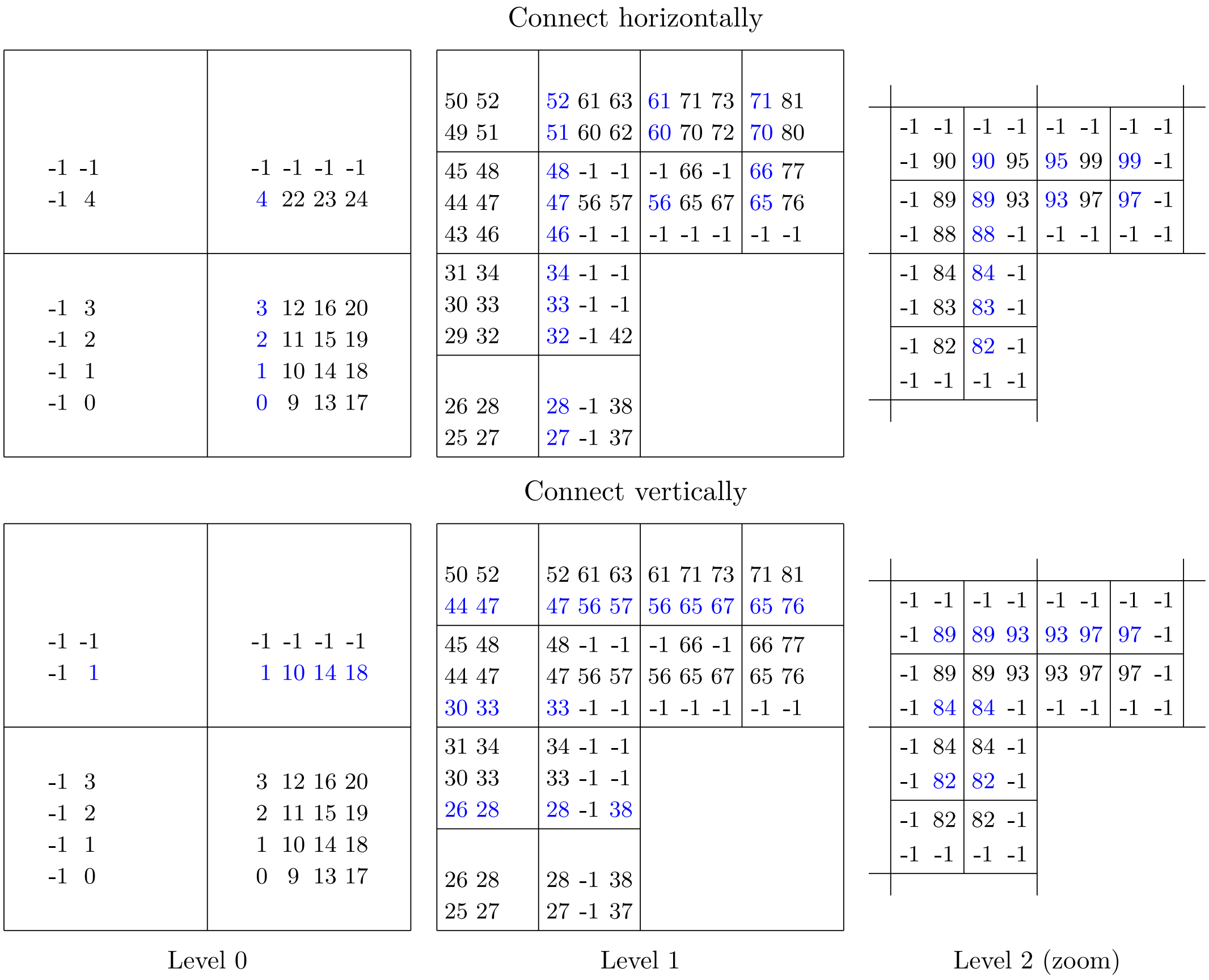

**Fig. 15.** Connect global ids over cell interfaces.

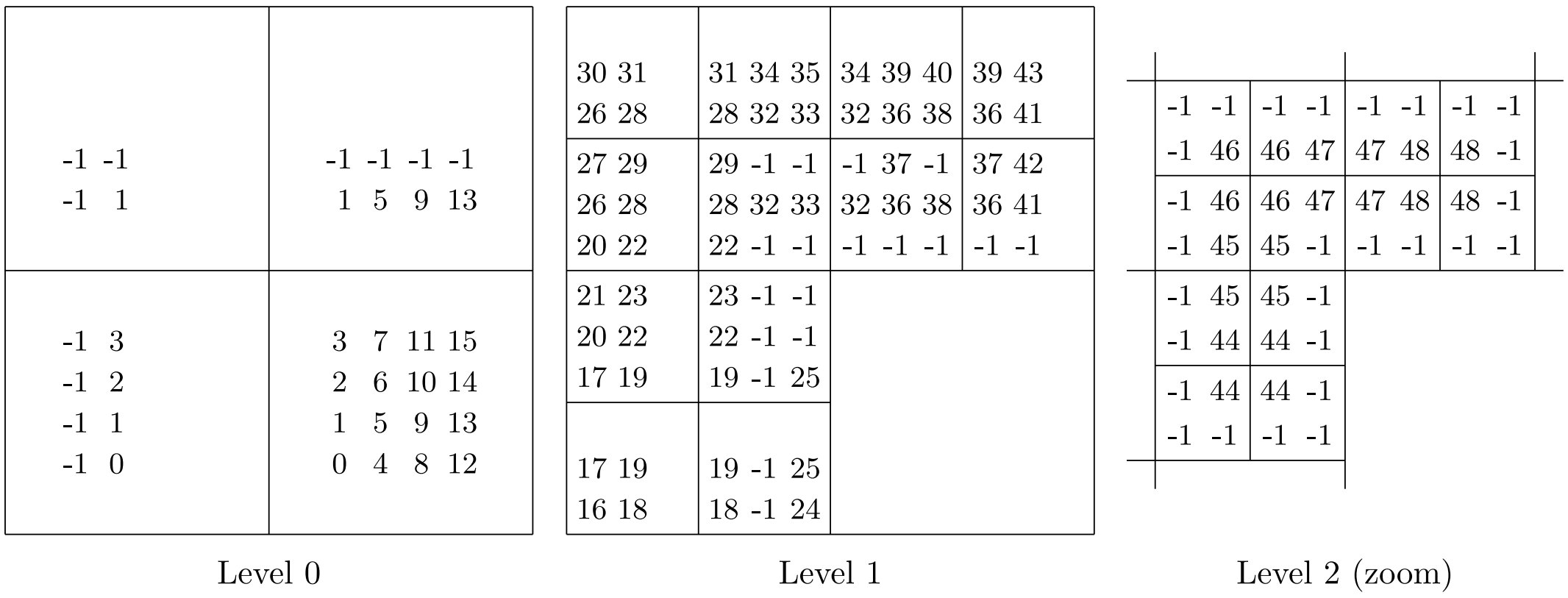

**Fig. 16.** Eliminate unassigned global basis function indices.

After the mesh and basis construction, we use the location maps of the multi-level *hp*-basis to allocate a sparse matrix in a compressed sparse row (CSR) format and a dense load vector. This ensures that our subsequent assembly

element 0 : [0, 1, 2, 3, 4, 5, 6, 7, 8, 9, 10, 11, 12, 13, 14, 15]
element 1 : [16, 17, 18, 19, 0, 1, 2, 3]
element 2 : [17, 20, 21, 19, 22, 23, 0, 1, 2, 3]
element 3 : [18, 19, 24, 25, 0, 1, 2, 3]
element 4 : [20, 26, 27, 22, 28, 29, 1]
element 5 : [26, 30, 28, 31, 1]
element 6 : [28, 31, 32, 34, 33, 35, 1]
element 7 : [32, 34, 36, 39, 38, 40, 1, 5, 9, 13]
element 8 : [36, 37, 41, 42, 1, 5, 9, 13]
element 9 : [36, 39, 41, 43, 1, 5, 9, 13]
element 10 : [44, 19, 22, 23, 25, 0, 1, 2, 3]
element 11 : [44, 45, 19, 22, 23, 25, 0, 1, 2, 3]
element 12 : [44, 19, 22, 23, 25, 0, 1, 2, 3]
element 13 : [44, 45, 19, 22, 23, 25, 0, 1, 2, 3]
element 14 : [45, 46, 22, 28, 29, 32, 33, 1]
element 15 : [46, 22, 28, 29, 32, 33, 1]
element 16 : [45, 46, 47, 22, 28, 29, 32, 33, 1]
element 17 : [46, 47, 22, 28, 29, 32, 33, 1]
element 18 : [47, 48, 32, 36, 37, 38, 1, 5, 9, 13]
element 19 : [47, 48, 32, 36, 37, 38, 1, 5, 9, 13]
element 20 : [48, 32, 36, 37, 38, 1, 5, 9, 13]
element 21 : [48, 32, 36, 37, 38, 1, 5, 9, 13]

**Fig. 17.** Element location maps including parent cells (colors represent parent levels). (For interpretation of the references to color in this figure legend, the reader is referred to the web version of this article.)

loop over all elements can add the element matrices to the global matrix without allocating new non-zero entries. For each element, we again obtain the current location map from the basis and allocate a dense element matrix and an element load vector according to the length of the location map. To determine the number of Gauss–Legendre points per direction, we query the basis for the maximum polynomial degrees of the element which are obtained from the sizes of the tensor-product masks of the element and its parents. For each of the integration points (usually $p+1$ per direction), the shape functions and their derivatives are evaluated and the finite element integrals are added to the element system. Finally, the element matrix and element load vector are assembled into the global system using the current location map.

We then use a diagonally preconditioned conjugate gradient method to solve the global linear system, which performs well for any number of refinement levels on our examples. In immersed settings, the multi-level $hp$-method has been combined with Additive-Schwarz [20] and multigrid preconditioners [21] (exploiting the hierarchical nature) or with sparse direct solvers [22]. After obtaining the solution coefficients, we write out VTU files for visualizing the results in Paraview. We capture the high-order nature of the approximation by evaluating the solution on a finer sample grid on each element that is chosen depending on the polynomial degree.

## 4. Numerical examples

In this section, we use our hierarchical framework to build several adaptive meshes. We implement the algorithms and data structures described in C++ for arbitrary dimensions. The source code for all examples is available at https://gitlab.com/hpfem/publications/2021_nd-mlhp.

### 4.1. Corner singularity

Our first example features a corner singularity in the origin that we refine towards. In two dimensions this benchmark is equivalent to computing one quadrant of an L-shaped domain and exploiting the symmetry of the solution (Fig. 18(a)). Similarly, the three-dimensional analog is computing one octant of a Fichera cube with a corner singularity (Fig. 18(b)). We define the following manufactured solution

$$u = \sqrt{r} \qquad \text{with} \qquad r = \sqrt{\sum_{i=0}^{d-1} x_i^2} \tag{6}$$

on $\Omega = [0, 1]^d$ with $d \geq 2$, as Fig. 19 shows. Substituting (6)(a) into (1) and using $\kappa = 1$ yields

$$f = \frac{3-2d}{4} r^{-3/2}.$$

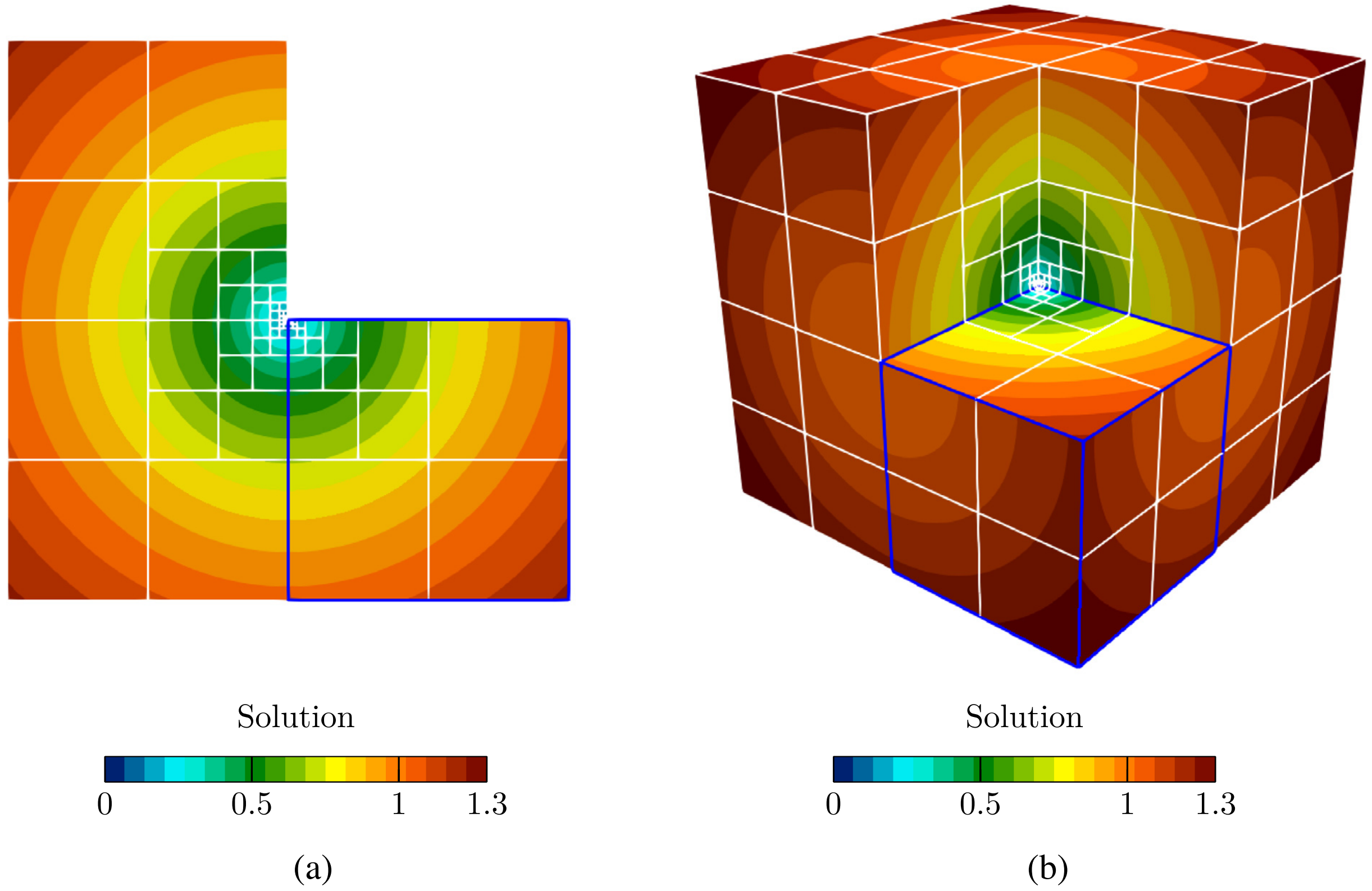

**Fig. 18.** (a) L-shaped domain benchmark. (b) Fichera corner benchmark.

We impose homogeneous Neumann boundary conditions on the sides intersecting the coordinate planes: the left and bottom edges in two dimensions, and the left, bottom, and front faces in three dimensions. On the remaining faces, we impose Dirichlet boundary conditions.

We discretize the domain $\Omega = [0, 1]^d$ using a base mesh with two elements in each direction. We perform two refinement studies where each new computation adds another refinement level towards the singularity at the origin. The first strategy uses a full tensor-product space with a uniform polynomial degree $p$ equal to $r + 1$, where $r$ is the maximum refinement level. The second strategy employs a trunk space with a linear grading from $p = r + 1$ on the base elements to $p = 1$ on the elements with refinement level $r$. Fig. 19(b) shows the graded polynomial degree distribution with one until four levels of refinement and the numerical solution using a uniform tensor-product space with $p = 5$ and $r = 4$. Fig. 20 shows the convergence results for different spatial dimensions.

While the convergence of the discretization error versus the number of unknowns is necessary to gain insight into the method's approximation quality, we ultimately believe the method's success depends almost exclusively on runtime performance and memory consumption. We now show that we achieve exponential convergence in terms of these crucial quantities. We use a diagonally preconditioned conjugate gradient method to solve the finite element system. Therefore, the number of non-zeros in the matrix dominates the memory consumption. The two most significant contributions to the runtime are the assembly and the solution of the linear system.

In the following, we analyze the asymptotic behavior of these quantities on the example. To allow comparisons with runtime measurements, we perform our analysis in terms of the refinement study index $r$. For brevity, we focus on the full tensor-product space using a uniform polynomial degree elevation. Each new refinement level adds $2^d$ cells that overlay the cell in the origin of the parent level. We can therefore compute the number of elements (leaf cells) as $r(2^d - 1) + 1 = \mathcal{O}(r)$.

Each leaf cell contributes to at most $(p + 1)^d \sim p^d$ basis functions defined on the same refinement level, but it also receives contributions from all coarser levels. The key to why we get good scaling is that parent cells only contribute basis functions active on their boundary. This results in at most $(p+1)^d - (p-1)^d \sim p^{d-1}$ functions per (coarser) level. Because $r \sim p$, we can conclude that these contributions scale with

$$\mathcal{O}\left(rp^{d-1}\right) \underset{r\sim p}{\longrightarrow} \mathcal{O}\left(p^d\right)$$

and hence are of the same order as the $(p-1)^d = \mathcal{O}\left(p^d\right)$ internal functions coming from the classical $p$-version of the finite element method.

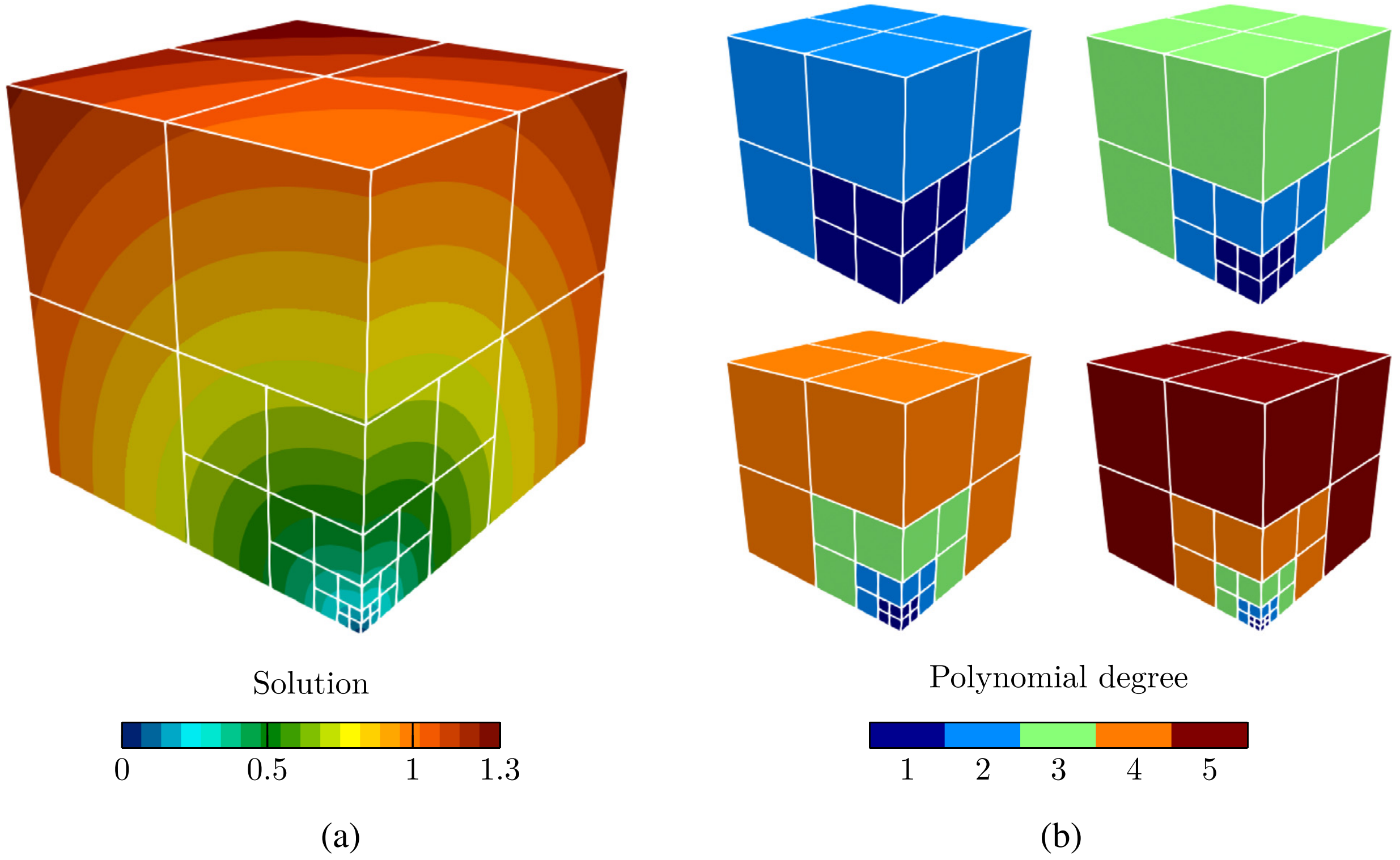

**Fig. 19.** (a) Singular benchmark solution with four refinement levels. (b) Graded polynomial degrees for refinement depths $r = 1 - 4$. (For interpretation of the references to color in this figure legend, the reader is referred to the web version of this article.)

We then bound the number of non-zeros in the assembled matrix to some value between the sum of all element matrix sizes (upper bound) and the sum of all internal contributions (lower bound). Using $\mathcal{O}\left(rp^{d-1} + p^d\right)$ as upper-bound for the number of element unknowns and multiplying with $O(r)$ number of elements yields

$$\mathcal{O}(r)\, \mathcal{O}\left(rp^{d-1} + p^d\right)^2 \underset{r\sim p}{\longrightarrow} \mathcal{O}\left(p^{2d+1}\right)$$

as an estimate of the complexity for the number of non-zeros in the assembled sparse matrix. The same result is obtained when considering only internal contributions.

Similarly, we estimate the assembly effort by summing up the number of operations needed to integrate each element matrix numerically. Assuming a standard Gauss–Legendre quadrature, we need $(p+1)^d = \mathcal{O}\left(p^d\right)$ integration points, on which we evaluate an outer product of $\mathcal{O}\left(rp^{d-1} + p^d\right)$ shape functions. Together, we obtain

$$\mathcal{O}(r)\, \mathcal{O}\left(p^d\right) \mathcal{O}\left(rp^{d-1} + p^d\right)^2 \underset{r\sim p}{\longrightarrow} \mathcal{O}\left(p^{3d+1}\right).$$

Fig. 21 shows that the number of non-zeros scales as the analysis above shows, while the number of conjugate gradient iterations increases roughly with $p^{d-1}$; the analysis of this empirical observation depends on the condition number of the sparse matrix that exceeds the scope of this work. We iterate until $\sqrt{r \cdot D(r)} < 10^{-10}$, where $D(r)$ is the diagonal preconditioner applied to the residual $r = Ax - b$. Moreover, we are not utilizing the symmetry of the sparse matrix, which significantly simplifies the parallel implementation of the sparse matrix–vector product. By combining the results for the number of non-zeros with the number of iterations, we estimate the runtime complexity for our iterative solver as $\mathcal{O}(p^{3d})$. Fig. 21 shows that our results verify our theoretical estimates. Although the assembly of the linear system scales with the highest order, it only starts significantly contributing to the total runtime towards the end of the curves. Moreover, we want to mention the possibility to construct direct solvers with a linear complexity with respect to the number of unknowns for the type of meshes used in this example [23].

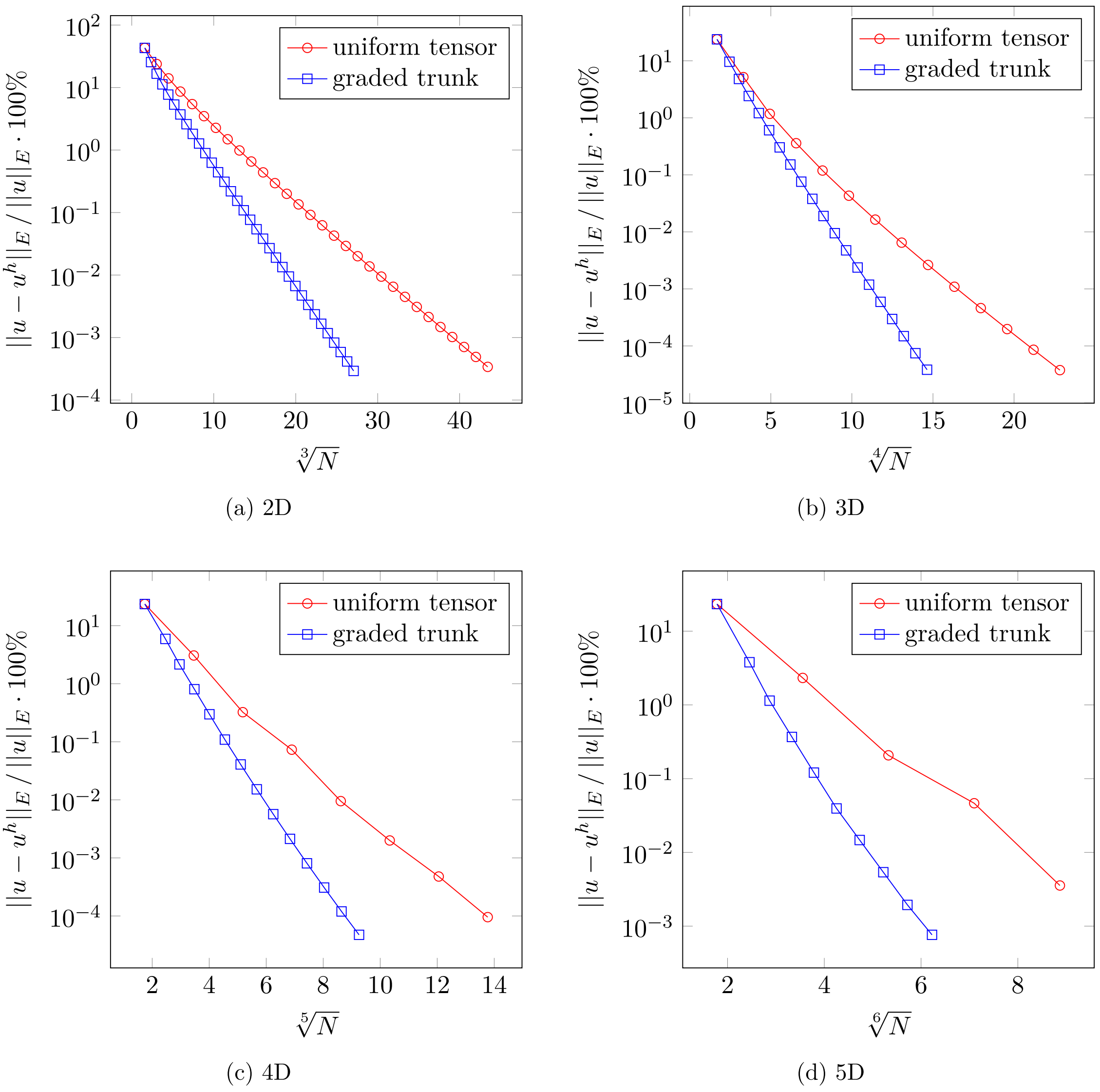

**Fig. 20.** Convergence of the error in the energy norm.

The measurements change significantly when considering the trunk space with the graded polynomial degree distribution. While the overall scaling is the same, the absolute effort per degree of freedom is higher, and the relative runtime of the assembly increases compared to the linear solution. We identify two aspects that are responsible for these differences. First, although the fine levels were assigned a lower polynomial degree, they still obtain contributions from the higher-order base elements. Due to these, even the finest level with $p = 1$ requires $r + 2$ Gauss–Legendre points per direction to integrate the parent contributions sufficiently accurate. Second, removing internal modes increases the connectivity in the trunk space, resulting in higher runtime and increased memory consumption per unknown. If we consider two interacting shape functions $N_i$ and $N_j$, then the effort for assembling this interaction is proportional to the number of elements $N_i$ and $N_j$ overlap. Internal functions overlap only on one element, but functions active on interfaces overlap at least two elements, making internal functions cheaper to assemble. After the assembly, however, this interaction results in a single entry in the sparse matrix regardless of the number of elements the supports overlap.

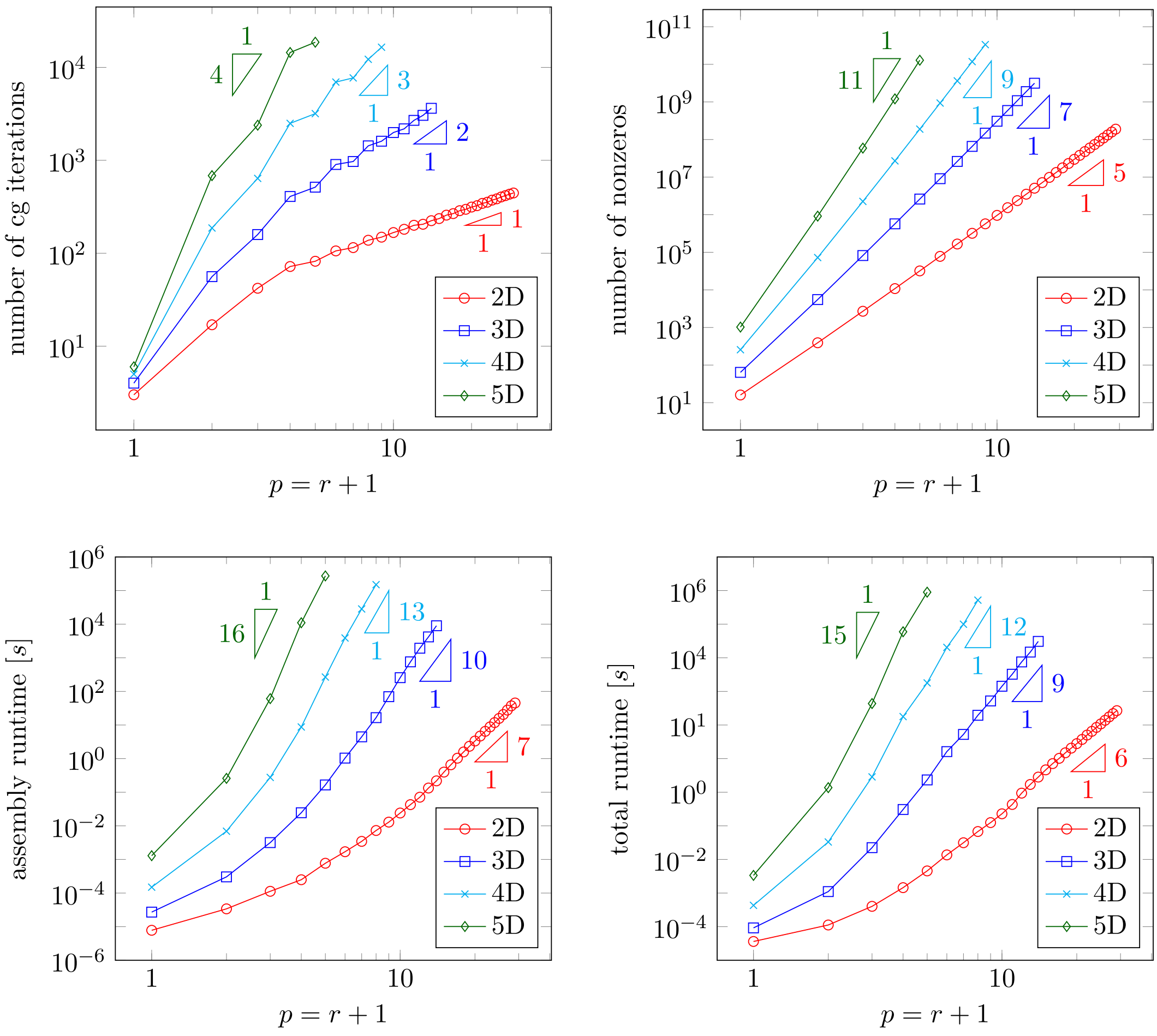

**Fig. 21.** Number of conjugate gradient iterations, number of non-zeros in the sparse matrix and serial runtime comparison for the uniformly graded tensor spaces. The total runtime excludes the error computation.

## 4.2. Transient problem with dynamic (de-)refinement

We demonstrate the effectiveness of the dynamic refinement and derefinement considering the transient version of (1) with constant coefficients $c$ and $\kappa$:

$$\begin{aligned} c\dot{u} - \kappa \Delta u &= f && \text{on } \Omega \times T \\ u &= u_0 && \text{at } t = 0 \\ u &= u_b && \text{on } \Gamma_D \\ \kappa \nabla u \cdot n &= 0 && \text{on } \Gamma_N. \end{aligned} \tag{7}$$

We discretize (7) using the method of Rothe, which discretizes time first and therefore allows for non-matching spatial discretizations in each time-step (as opposed to the method of lines, which discretizes space first). We rewrite the domain part of (7) as $c\dot{u} = f + \kappa \Delta u$ and apply the general-$\theta$ method:

$$\frac{c}{\Delta t}\left(u^{n+1} - u^n\right) = \theta\left(f^{n+1} + \kappa \Delta u^{n+1}\right) + (1-\theta)\left(f^n + \kappa \Delta u^n\right), \tag{8}$$

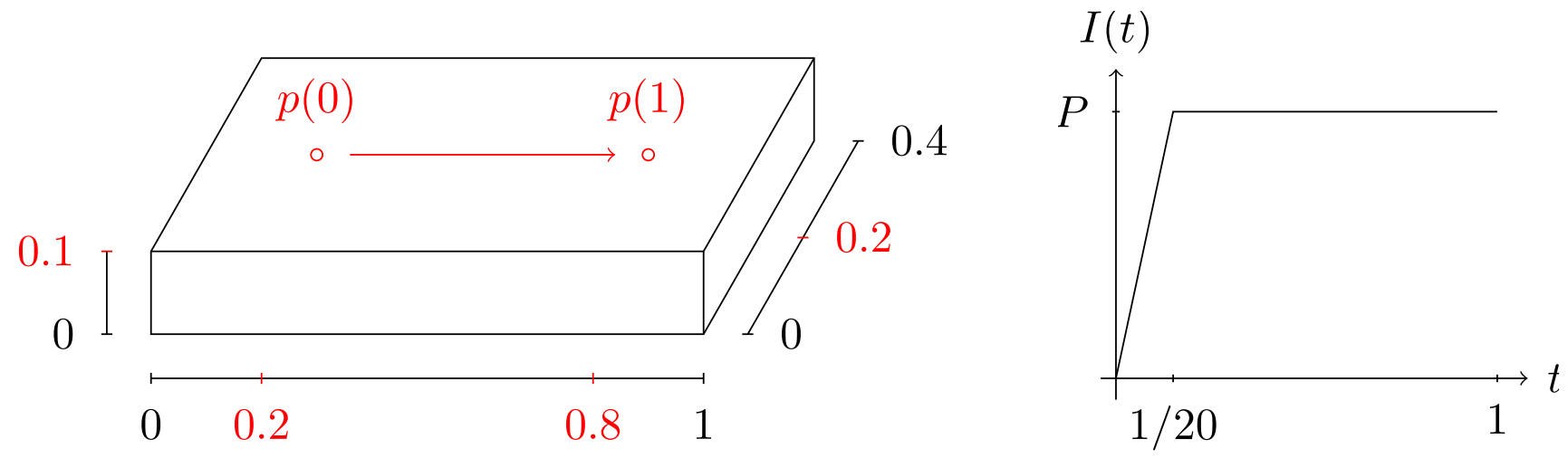

| | |
|---|---|
| $c$ | 1 |
| $\kappa$ | 0.01 |
| $u_0$ | 25 |
| $\sigma$ | 0.0021 |
| $P$ | 0.968 |

**Fig. 22.** Domain and path definitions, intensity function and problem parameters.

where $\theta = 0$ and $\theta = 1$ yield the first-order explicit and implicit Euler methods and $\theta = 1/2$ yields the second-order implicit midpoint (Crank–Nicolson) rule. Reorganizing the terms in (8) yields

$$\frac{c}{\Delta t} u^{n+1} - \kappa\, \theta\, \Delta u^{n+1} = \frac{c}{\Delta t} u^n + \kappa\, (1-\theta)\, \Delta u^n + \theta f^{n+1} + (1-\theta) f^n. \tag{9}$$

We then discretize the weak form of (9) to obtain a system of equations $K_{ij}^{n+1}\, \hat{u}_j^{n+1} = F_i^{n+1}$, where

$$\begin{aligned} K_{ij}^{n+1} &= \int_\Omega \frac{c}{\Delta t} N_i N_j + \kappa\, \theta\, \nabla N_i \cdot \nabla N_j \,\mathrm{d}\Omega \\ F_i^{n+1} &= \int_\Omega N_i \left( \frac{c}{\Delta t} u^n + \theta f^{n+1} + (1-\theta) f^n \right) - \kappa\, (1-\theta)\, \nabla N_i \cdot \nabla u^n \,\mathrm{d}\Omega. \end{aligned} \tag{10}$$

The basis functions $N_i$ and $N_j$ in (10) belong to the time-slice $t^{n+1}$. The old solution $u^n$ and its derivatives $\nabla u^n$, however, belong to the time-slice $t^n$. The representation of the solution in two time-slices requires our implementation to provide two bases (for $t^n$ and for $t^{n+1}$) to assemble the equation system resulting from one time-step. Furthermore, we ensure $C^\infty$ continuity on each integration cell by distributing Gauss–Legendre points on the intersection of both grids.

Now, we demonstrate the performance and consistency of the two-time-slice representation on an example motivated by the simulation of powder bed fusion processes in additive manufacturing, where a laser selectively melts the metal powder in a layer-by-layer fashion. We choose a volumetric heat source in the shape of a Gaussian bell function that travels along a path $p(t)$ with varying intensity $I(t)$. We therefore have

$$f(x,t) = I(t)\, q(x - p(t)),$$

with

$$q(x) = \frac{1}{(2\pi)^{d/2}\, \sigma^d} \exp\left( \frac{-x_i^2}{2\sigma^2} \right).$$

Assuming that $I(t) = 0$ for $t < 0$ and $u(x,0) = u_0 = \text{const}$, we obtain the following semi-analytical solution on an infinite domain:

$$\begin{aligned} u(x,t) &= \int_0^t I(\tau)\, g(x - p(\tau),\, t - \tau)\, \mathrm{d}\tau + u_0 \\ g(x,t) &= \frac{1}{c\,(2\pi)^{d/2}\, w(t)^d} \exp\left( \frac{-x_i^2}{2 w(t)^2} \right) \\ w(t) &= \sqrt{\sigma^2 + \frac{2\kappa t}{c}}. \end{aligned} \tag{11}$$

Here, $g$ represents the solution to (7) with $f = q(x)\, \delta(t)$. We then restrict $\Omega$ and $T$ to a three-dimensional cube on the time interval [0, 1] and select a linear source path on the top face. Fig. 22 shows the geometric setup and the physical parameters we use. To conform to (11), we choose homogeneous boundary conditions on the top face due to symmetry and Dirichlet boundary conditions on all other faces.

Fig. 23 shows the numerical solution of this setup after 232 of the total 256 Crank–Nicolson steps. We refine the mesh of each time-step five levels towards the current source center and then gradually coarsen it along the previous source positions. The construction of this refinement can be found in the source repository and is explained for general laser paths in additive manufacturing in [12]. We assign the polynomial degrees 5, 5, 4, 3, 3, 2 to the refinement levels 0 to 5, resulting in an average of 3850 unknowns per time-step. The computation was performed on a Desktop computer using an i9-9900KF processor with all eight cores running at 4.6 GHz. We configured our code to use 32-bit unsigned integers for cell indices and sparse (= degree of freedom) indices and we used the g++-10 compiler.

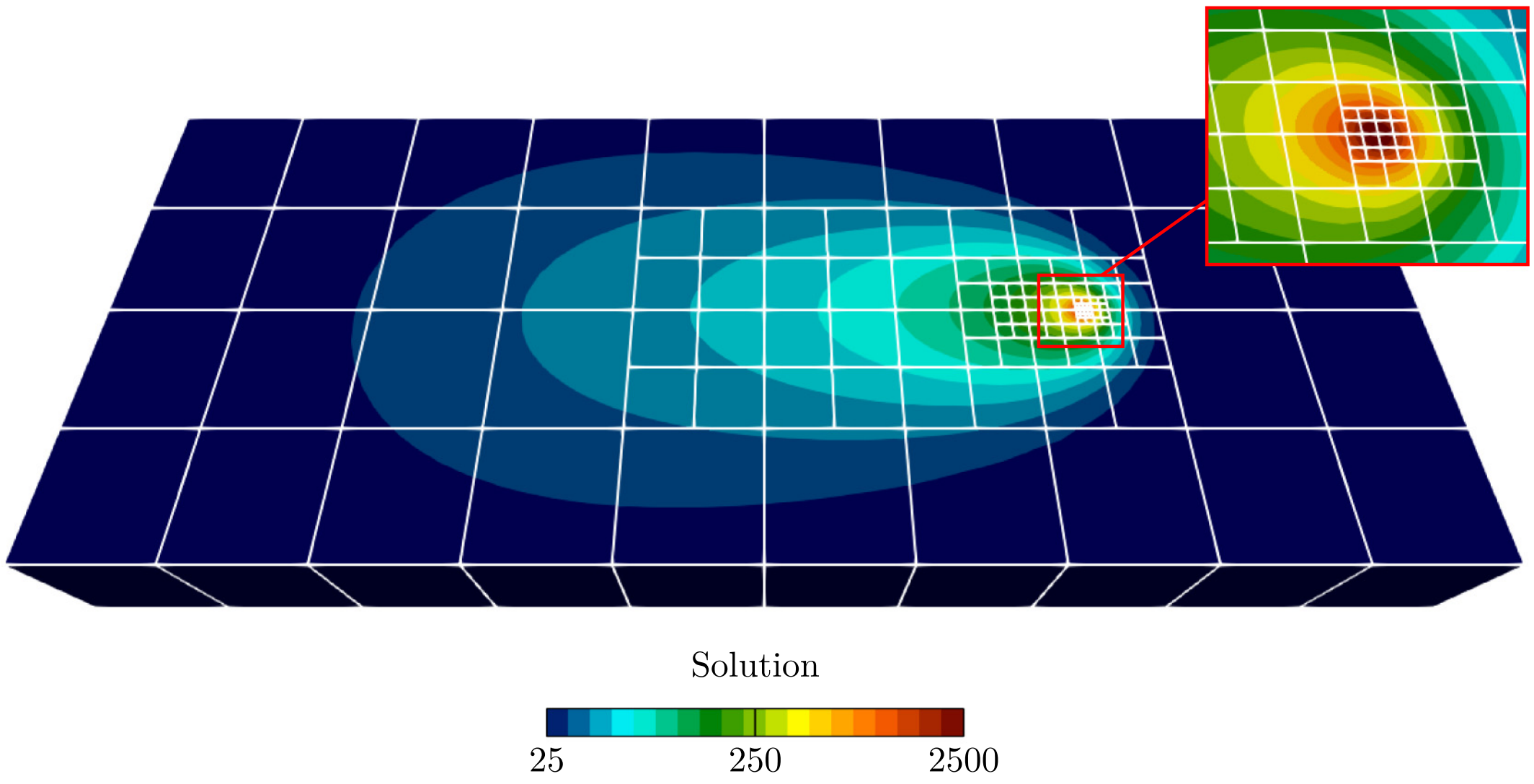

**Fig. 23.** Solution after 232 time-steps using a logarithmic color map. (For interpretation of the references to color in this figure legend, the reader is referred to the web version of this article.)

The computation took 146.9 s for a single-threaded execution and 23.6 s in a shared-memory environment using eight threads, resulting in a speedup of 6.2x. In the parallel execution, a single time-step computes on average in 92 ms, which splits into 1 ms for creating and refining the mesh and constructing the multi-level $hp$-basis, 50 ms for allocating the sparse matrix and assembling the linear system, and 22 ms for the iterative linear solution using a conjugate-gradient method. An additional 18 ms is spent on computing the Dirichlet boundary condition on the left, right, front, back and bottom faces according to (11) using a global $L^2$ projection on the boundary. This is required to evaluate the discretization error, although integrating the semi-analytical solution on each quadrature point is expensive. The memory required for storing the sparse matrix averaged around 6 MB, while the topological information and the basis data together used around 75 kB of memory. From these numbers we conclude that the presented algorithms and data structures perform well in runtime (1 ms out of 92 ms) and memory consumption (compared to the sparse matrix) and that there is almost no sequential execution and only little communication overhead for practical problem sizes. We also believe that our multi-level $hp$-framework is well-suited for massively parallel computations in distributed memory environments, which we plan to address in future work.

We then perform convergence studies for spatial and temporal discretizations in terms of the error in the $L^2$ norm. In Fig. 24(a), we increase the polynomial degree for a fixed temporal discretization of 4096 Crank–Nicolson steps. $\tilde{N}$ denotes the average number of unknowns per time-step. Fig. 24(b) shows the convergence of different time marching schemes for a fixed spatial discretization using a trunk space with $p = 12$. For all cases, we use the same locally refined meshes. We evaluate the norms numerically using a trapezoidal rule in time and a Gauss–Legendre quadrature in space with $p + 1$ points per direction on each integration cell. The integral in the semi-analytical solution is again evaluated numerically on each integration point in a stepwise fashion with $dt \approx 1/100$ and 30 Gauss–Legendre points per interval.

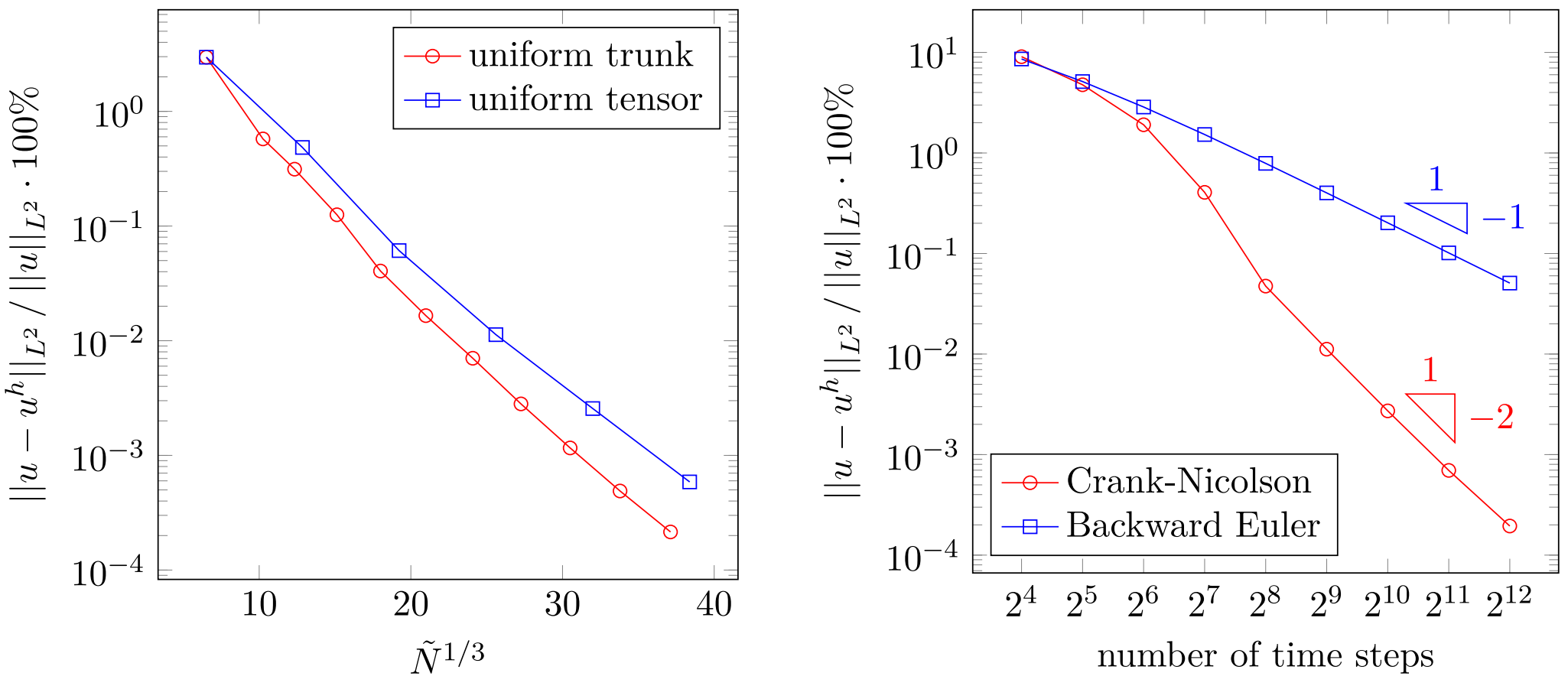

(a) Spatial $p$-refinement. (b) Temporal refinement.

**Fig. 24.** Convergence of spatial and temporal discretizations.

## 5. Conclusion

We describe a data-oriented construction for multi-level $hp$-bases. The algorithmic framework only operates on the $d - 1$ dimensional slices of the element data, which significantly simplifies the implementation of the method for higher dimensions. While we focus on Cartesian grids here, these data structures and algorithms easily generalize to unstructured meshes. We consider a Fichera corner example to show that our implementation scales well, even for small problem sizes. However, the exponential complexity in the number of spatial dimensions makes high-dimensional computations very expensive regarding runtime and memory consumption. Finally, we apply our methodology to a transient heat equation and show the potential of dynamic local refinement and derefinement in bridging the vastly different spatial scales arising, for example, in powder bed fusion processes. Our next step will be to address spatial and temporal scales using the presented method combined with space–time finite element methods.

## Declaration of competing interest

The authors declare the following financial interests/personal relationships which may be considered as potential competing interests: Philipp Kopp reports financial support that was provided by Deutsche Forschungsgemeinschaft (DFG). Victor M. Calo reports financial support that was provided by CSIRO Professorial Chair in Computational Geoscience at Curtin University. Victor M. Calo reports was provided by European Union's Horizon 2020 research and innovation programme.

## Data availability

We provide a link to our source code repository that includes all our examples.

## Acknowledgments

This research was funded by the Deutsche Forschungsgemeinschaft (DFG, German Research Foundation) – grant KO 4570/2-1. We greatfully aknowledge this support. This publication was also made possible in part by the Professorial Chair in Computational Geoscience at Curtin University and the Deep Earth Imaging Enterprise Future Science Platforms of the Commonwealth Scientific Industrial Research Organisation, CSIRO, of Australia. This project has received funding from the European Union's Horizon 2020 research and innovation programme under the Marie Sklodowska-Curie grant agreement No 777778 (MATHROCKS). The Curtin Corrosion Centre and the Curtin Institute for Computation kindly provide ongoing support.

## Appendix. Construction of initial tensor-product masks for the trunk space

Section 3 mentions that using all functions from the tensor-product is not always necessary. If we consider a monomial basis

$$\left\{ r_0^{\alpha_0} \cdot \ldots \cdot r_{d-1}^{\alpha_{d-1}} \;\middle|\; \text{for } \alpha_i = 0, \ldots, p_i \right\}, \tag{12}$$

then picking only functions where $\sum \alpha_i \leq \max p_i$ will not reduce the convergence order. Fig. 25 shows the monomials for a tensor-product space with $p = (3, 4)$ compared to the reduced basis containing only polynomials with a combined order smaller or equal to $\max p_i = 4$.

As monomials are not suitable for constructing compatible finite element bases, we transfer this concept to integrated Legendre functions. In Section 3, we started our algorithm with a full tensor-product mask for leaf cells. To construct initial trunk space masks, we only activate entries of $M_{\boldsymbol{\alpha}}$ if $\sum \alpha_i \leq \max p_i$. The integrated Legendre functions are hierarchical except in the first two functions ($I_0$ and $I_1$ are linear). As a result, we need to ensure that in no case a function with component $I_0(r_i)$ is active, while the corresponding function with component $I_1(r_i)$ is inactive. Therefore, we take the first slice for each coordinate axis and copy its entries to the second slice. Figs. 26 and 27 demonstrate this procedure on two- and three-dimensional examples.

While using a trunk space does not reduce the order of convergence, it can significantly reduce the accuracy for a given setup. When increasing the polynomial degree to compensate, one generally obtains fewer degrees of freedom for a given accuracy, at the cost of having higher connectivity per unknown. With higher dimensions, the difference between both becomes more pronounced as the limit for $p \to \infty$ yields $d!$ times fewer unknowns per element (assuming uniform $p$ in all directions). This limit derives from considering an $n$-simplex in $d$ dimensional

tensor product space

| | | | |
|---|---|---|---|
| $r_1^4$ | $r_0 r_1^4$ | $r_0^2 r_1^4$ | $r_0^3 r_1^4$ |
| $r_1^3$ | $r_0 r_1^3$ | $r_0^2 r_1^3$ | $r_0^3 r_1^3$ |
| $r_1^2$ | $r_0 r_1^2$ | $r_0^2 r_1^2$ | $r_0^3 r_1^2$ |
| $r_1$ | $r_0 r_1$ | $r_0^2 r_1$ | $r_0^3 r_1$ |
| $1$ | $r_0$ | $r_0^2$ | $r_0^3$ |

monomial trunk space

| | | | |
|---|---|---|---|
| $r_1^4$ | | | |
| $r_1^3$ | $r_0 r_1^3$ | | |
| $r_1^2$ | $r_0 r_1^2$ | $r_0^2 r_1^2$ | |
| $r_1$ | $r_0 r_1$ | $r_0^2 r_1$ | $r_0^3 r_1$ |
| $1$ | $r_0$ | $r_0^2$ | $r_0^3$ |

**Fig. 25.** Monomials for $p = (3, 4)$.

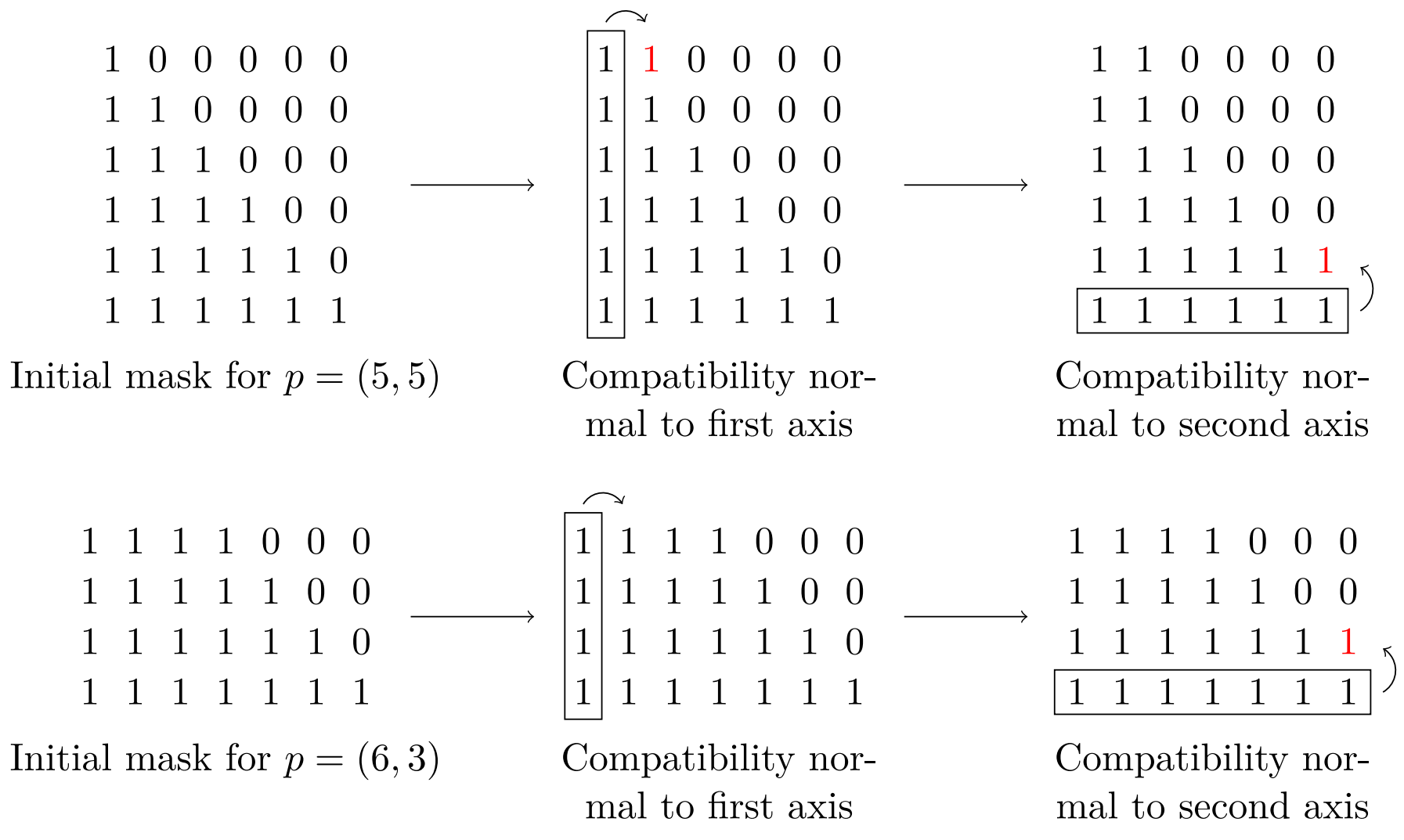

**Fig. 26.** Construction of two-dimensional trunk space initial tensor-product masks.

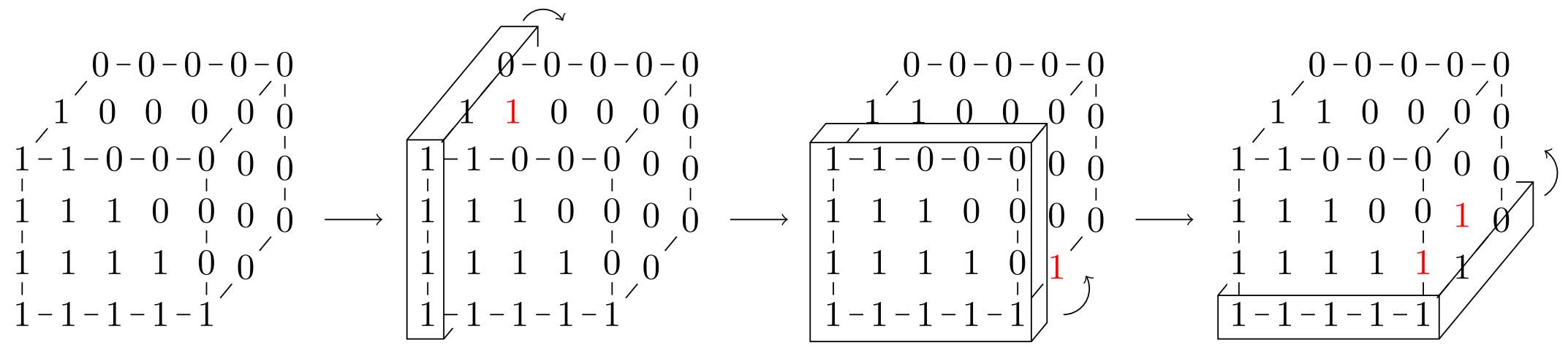

**Fig. 27.** Construction of three-dimensional trunk space initial tensor-product masks.

space with one vertex at the origin and one vertex on each coordinate axis with distance 1. The volume of this simplex is equal to $1/d!$, while the volume of the corresponding unit-$n$-cube is 1.